\documentclass[a4paper,11pt,draft]{amsart}
\usepackage{amssymb,amscd}
\usepackage[cp1251]{inputenc}
\usepackage[T2A]{fontenc}
\newtheorem{theo}{Theorem}[section]
\newtheorem{prop}[theo]{Proposition}
\newtheorem{coro}[theo]{Corollary}
\newtheorem{lemm}[theo]{Lemma}

\theoremstyle{definition}

\newtheorem{exam}[theo]{Example}

\newtheorem{cons}[theo]{Construction}

\theoremstyle{remark}

\numberwithin{equation}{section}

\newcommand{\mb}[1]{{\textbf {\textit#1}}}

\def\C{\mathbb C}
\newcommand{\R}{\mathbb R}
\newcommand{\Q}{\mathbb Q}
\newcommand{\Z}{\mathbb Z}
\renewcommand{\k}{\mathbf k}

\def\D{\Delta}
\def\L{\Lambda}
\def\b{\beta}
\def\e{\varepsilon}
\def\phi{\varphi}
\def\l{\lambda}
\def\s{\sigma}

\def\zk{\mathcal Z_K}

\newcommand{\bideg}{\mathop{\rm bideg}}
\newcommand{\cc}{\mathop{\rm cc}}
\newcommand{\cone}{\mathop{\rm cone}}
\newcommand{\djs}{\mathop{\mbox{\it DJ\/}}\nolimits}
\newcommand{\link}{\mathop{\rm link}\nolimits}
\renewcommand{\star}{\mathop{\rm star}\nolimits}
\def\Tor{\operatorname{Tor}}
\def\trk{\operatorname{trk}}

\def\id{\mathrm{id}}

\begin{document}

\title{Cohomology of face rings, and torus actions}

\author{Taras Panov}\thanks{The author was supported by an LMS grant for young Russian mathematicians
at the University of Manchester, and also by the Russian
Foundation for Basic Research, grant no.~04-01-00702.}
\address{Department of Geometry and Topology, Faculty of Mathematics and Mechanics, Moscow State University,
Leninskiye Gory, Moscow 119992, Russia\newline \emph{and
}Institute for Theoretical and Experimental Physics, Moscow
117259, Russia}
\email{tpanov@mech.math.msu.su}

\begin{abstract}
In this survey article we present several new developments of
`toric topology' concerning the cohomology of face rings (also
known as Stanley--Reisner algebras). We prove that the integral
cohomology algebra of the \emph{moment-angle complex} $\zk$
(equivalently, of the complement $U(K)$ of the coordinate subspace
arrangement) determined by a simplicial complex $K$ is isomorphic
to the $\Tor$-algebra of the face ring of~$K$. Then we analyse
Massey products and formality of this algebra by using a
generalisation of Hochster's theorem. We also review several
related combinatorial results and problems.
\end{abstract}

\maketitle

\section{Introduction}
This article centres on the cohomological aspects of `toric
topology', a new and actively developing field on the borders of
equivariant topology, combinatorial geometry and commutative
algebra. The algebro-geometric counterpart of toric topology,
known as `toric geometry' or algebraic geometry of \emph{toric
varieties}, is now a well established field in algebraic geometry,
which is characterised by its strong links with combinatorial and
convex geometry (see the classical survey paper~\cite{dani78} or
more modern exposition~\cite{fult93}). Since the appearance of
Davis and Januszkiewicz's work~\cite{da-ja91}, where the concept
of a \emph{(quasi)toric manifold} was introduced as a topological
generalisation of smooth compact toric variety, there has grown an
understanding that most phenomena of smooth toric geometry may be
modelled in the purely topological situation of smooth manifolds
with a nicely behaved torus action.

One of the main results of~\cite{da-ja91} is that the equivariant
cohomology of a toric manifold can be identified with the
\emph{face ring} of the quotient simple polytope, or, for more
general classes of torus actions, with the face ring of a certain
simplicial complex~$K$. The ordinary cohomology of a quasitoric
manifold can also be effectively identified as the quotient of the
face ring by a \emph{regular sequence} of degree-two elements,
which provides a generalisation to the well-known
Danilov--Jurkiewicz theorem of toric geometry. The notion of the
face ring of a simplicial complex sits in the heart of Stanley's
`Combinatorial commutative algebra'~\cite{stan96}, linking
geometrical and combinatorial problems concerning simplicial
complexes with commutative and homological algebra. Our concept of
toric topology aims at extending these links and developing new
applications by applying the full strength of the apparatus of
equivariant topology of torus actions.

The article surveys certain new developments of toric topology
related to the cohomology of face rings. Introductory remarks can
be found at the beginning of each section and most subsections. A
more detailed description of the history of the subject, together
with an extensive bibliography, can be found in~\cite{bu-pa02} and
its extended Russian version~\cite{bu-pa04}.

The current article represents the work of the algebraic topology
and combinatorics group at the Department of Geometry and
Topology, Moscow State University, and the author thanks all its
members for the collaboration and insight gained from numerous
discussions, particularly mentioning Victor Buchstaber, Ilia
Baskakov, and Arseny Gadzhikurbanov. The author is also grateful
to Nigel Ray for several valuable comments and suggestions that
greatly improved this text and his hospitality during the visit to
Manchester sponsored by an LMS grant.

\section{Simplicial complexes and face rings}\label{fring}
The notion of the face ring $\k[K]$ of a simplicial complex $K$ is
central to the algebraic study of triangulations. In this section
we review its main properties, emphasising functoriality with
respect to simplicial maps. Then we introduce the bigraded
$\Tor$-algebra $\Tor_{\k[v_1,\dots,v_m]}(\k[K],\k)$ through a
finite free resolution of $\k[K]$ as a module over the polynomial
ring. The corresponding bigraded Betti numbers are important
combinatorial invariants of~$K$.

\subsection{Definition and main properties}
Let $K=K^{n-1}$ be an arbitrary $(n-1)$-dimensional simplicial
complex on an $m$-element vertex set $V$, which we usually
identify with the set of ordinals $[m]=\{1,\ldots,m\}$. Those
subsets $\sigma\subseteq V$ belonging to $K$ are referred to as
\emph{simplices}; we also use the notation $\sigma\in K$. We count
the empty set $\varnothing$ as a simplex of~$K$. When it is
necessary to distinguish between combinatorial and geometrical
objects, we denote by $|K|$ a \emph{geometrical realisation} of
$K$, which is a triangulated topological space.

Choose a ground commutative ring $\k$ with unit (we are mostly
interested in the cases $\k=\Z,\Q$ or finite field). Let
$\k[v_1,\ldots,v_m]$ be the graded polynomial algebra over $\k$
with $\deg v_i=2$.  For an arbitrary subset
$\omega=\{i_1,\dots,i_k\}\subseteq[m]$, denote by $v_ \omega$ the
square-free monomial $v_{i_1}\dots v_{i_k}$.

The \emph{face ring} (or
\emph{Stanley--Reisner algebra}) of $K$ is the quotient ring
\[
  \k[K]=\k[v_1,\ldots,v_m]/\mathcal I_K,
\]
where $\mathcal I_K$ is the homogeneous ideal generated by all
monomials $v_\sigma$ such that $\sigma$ is not a simplex of $K$.
The ideal $\mathcal I_K$ is called the \emph{Stanley--Reisner
ideal} of~$K$.

\begin{exam}
Let $K$ be a 2-dimensional simplicial complex shown on
Figure~\ref{figfr}. Then
\[
  \k[K]=\k[v_1,\ldots,v_5]/(v_1v_5,v_3v_4,v_1v_2v_3,v_2v_4v_5).
\]
\begin{figure}[h]
\begin{center}
\begin{picture}(60,35)
\put(5,5){\line(1,0){50}} \put(5,5){\line(1,1){25}}
\put(5,5){\line(2,1){20}} \put(30,30){\line(1,-1){25}}
\put(25,15){\line(1,0){10}} \put(25,15){\line(1,3){5}}
\put(35,15){\line(-1,3){5}} \put(35,15){\line(2,-1){20}}
\put(25,25){\line(0,-1){10}} \put(20,20){\line(0,-1){7.5}}
\put(15,15){\line(0,-1){5}} \put(10,10){\line(0,-1){2.5}}
\put(7.5,7.5){\line(0,-1){1.25}}
\put(12.5,12.5){\line(0,-1){3.75}}
\put(17.5,17.5){\line(0,-1){6.25}}
\put(22.5,22.5){\line(0,-1){8.75}} \put(27.5,27.5){\line(0,-1){5}}
\put(35,25){\line(0,-1){10}} \put(40,20){\line(0,-1){7.5}}
\put(45,15){\line(0,-1){5}} \put(50,10){\line(0,-1){2.5}}
\put(52.5,7.5){\line(0,-1){1.25}}
\put(47.5,12.5){\line(0,-1){3.75}}
\put(42.5,17.5){\line(0,-1){6.25}}
\put(37.5,22.5){\line(0,-1){8.75}} \put(32.5,27.5){\line(0,-1){5}}
\put(5,5){\circle*{1}} \put(25,15){\circle*{1}}
\put(35,15){\circle*{1}} \put(30,30){\circle*{1}}
\put(55,5){\circle*{1}} \put(3,1){1} \put(56,1){3} \put(29,31){2}
\put(24,11){4} \put(34,11){5}
\end{picture}
\caption{} \label{figfr}
\end{center}
\end{figure}
\end{exam}

Despite its simple construction, the face ring appears to be a
very powerful tool allowing
us 
to translate the combinatorial properties of different particular
classes of simplicial complexes into the language of commutative
algebra. The resulting field of `Combinatorial commutative
algebra', whose foundations were laid by Stanley in his
monograph~\cite{stan96}, has attracted a lot of interest from both
combinatorialists and commutative algebraists.

Let $K_1$ and $K_2$ be two simplicial complexes on the vertex sets
$[m_1]$ and $[m_2]$ respectively. A set map $\phi\colon[m_1]\to
[m_2]$ is called a \emph{simplicial map} between $K_1$ and $K_2$
if $\phi(\sigma)\in K_2$ for any $\sigma\in K_1$; we often
identify such $\phi$ with its restriction to $K_1$ (regarded as a
collection of subsets of $[m_1]$), and use the notation
$\phi\colon K_1\to K_2$.

\begin{prop}\label{frmap}
Let $\phi\colon K_1\to K_2$ be a simplicial map. Define a map
$\phi^*\colon \k[w_1,\ldots,w_{m_2}]\to\k[v_1,\ldots,v_{m_1}]$ by
\[
  \phi^*(w_j):=\sum_{i\in\phi^{-1}(j)}v_i.
\]
Then $\phi^*$ induces a homomorphism $\k[K_2]\to\k[K_1]$, which we
will also denote by $\phi^*$.
\end{prop}
\begin{proof}
We have to check that $\phi^*(\mathcal I_{K_2})\subseteq\mathcal
I_{K_1}$. Suppose $\tau=\{j_1,\ldots,j_s\}\subseteq[m_2]$ is not a
simplex of~$K_2$. Then
\begin{equation}
\label{fstsum}
  \phi^*(w_{j_1}\cdots w_{j_s})=\sum_{i_1\in \phi^{-1}(j_1),\ldots,
  i_s\in \phi^{-1}(j_s)}v_{i_1}\cdots v_{i_s}.
\end{equation}
We claim that $\s=\{i_1,\ldots,i_s\}$ is not a simplex of $K_1$
for any monomial $v_{i_1}\cdots v_{i_s}$ in the right hand side of
the above identity. Indeed, if $\sigma\in K_1$, then
$\phi(\s)=\tau\in K_2$ by the definition of simplicial map, which
leads to a contradiction. Hence, the right hand side
of~\eqref{fstsum} is in~$\mathcal I_{K_1}$.
\end{proof}

\subsection{Cohen--Macaulay rings and complexes}
{\sloppy Cohen--Macaulay rings and modules play an important role
in homological commutative algebra and algebraic geometry. A
standard reference for the subject is~\cite{br-he98}, where the
reader may find proofs of the basic facts about Cohen--Macaulay
rings and regular sequences mentioned in this subsection. In the
case of simplicial complexes, the Cohen--Macaulay property of the
corresponding face rings leads to important combinatorial and
topological consequences.

}
Let $A=\oplus_{i\ge0}A^i$ be a finitely-generated commutative
graded algebra over~$\k$. We assume that $A$ is connected
($A^0=\k$) and has only even-degree graded components, so that we
do not need to distinguish between graded and non-graded
commutativity. We denote by $A_+$ the positive-degree part of $A$
and by $\mathcal H(A_+)$ the set of homogeneous elements in $A_+$.

A sequence $t_1,\ldots,t_n$ of algebraically independent
homogeneous elements of $A$ is called an \emph{hsop} (homogeneous
system of parameters) if $A$ is a finitely-generated
$\k[t_1,\ldots,t_n]$-module (equivalently, $A/(t_1,\ldots,t_n)$
has finite dimension as a $\k$-vector space).

\begin{lemm}[N\"other normalisation lemma]\label{noether}
Any finitely-generated graded algebra $A$ over a field $\k$ admits
an hsop. If $\k$ has characteristic zero and $A$ is generated by
degree-two elements, then a degree-two hsop can be chosen.
\end{lemm}

A degree-two hsop is called an \emph{lsop} (linear system of
parameters).

A sequence $\mb t=t_1,\ldots,t_k$ of elements of $\mathcal H(A_+)$
is called a \emph{regular sequence} if $t_{i+1}$ is not a zero
divisor in $A/(t_1,\ldots,t_i)$ for $0\le i<k$. A regular sequence
consists of algebraically independent elements, so it generates a
polynomial subring in $A$. It can be shown that $\mb t$ is a
regular sequence if and only if $A$ is a \emph{free}
$\k[t_1,\ldots,t_k]$-module.

An algebra 
$A$ is called \emph{Cohen--Macaulay} if it admits a
regular hsop $\mb t$. It follows that $A$ is Cohen--Macaulay if
and only if it is
a 
free
and 
finitely generated module over its polynomial subring.  If $\k$ is a
field of zero characteristic and $A$ is generated by degree-two
elements, then one can choose $\mb t$ to be an lsop. A simplicial
complex $K$ is called \emph{Cohen--Macaulay} (over $\k$) if its face
ring $\k[K]$ is Cohen--Macaulay.

\begin{exam}
Let $K=\partial\D^2$ be the boundary of a 2-simplex. Then
\[
  \k[K]=\k[v_1,v_2,v_3]/(v_1v_2v_3).
\]
The elements $v_1,v_2\in\k[K]$ are algebraically independent,
but do not form an hsop, since $\k[K]/(v_1,v_2)\cong\k[v_3]$ is
not finite-dimensional as a $\k$-space.  On the other hand, the
elements $t_1=v_1-v_3$, $t_2=v_2-v_3$ of $\k[K]$ form an hsop,
since $\k[K]/(t_1,t_2)\cong\k[t]/t^3$. It is easy to see that
$\k[K]$ is a free $\k[t_1,t_2]$-module with one 0-dimensional
generator~1, one 1-dimensional generator~$v_1$, and one
2-dimensional generator~$v_1^2$. Thus, $\k[K]$ is Cohen--Macaulay
and $(t_1,t_2)$ is a regular sequence.
\end{exam}

For
an 
arbitrary simplex $\s\in K$ define its \emph{link} and
\emph{star} as the subcomplexes
\begin{align*}
  \link_K\s&=\{\tau\in K\colon \s\cup\tau\in K,\;\s\cap\tau=\varnothing
  \};\\
  \star_K\s&=\{\tau\in K\colon \s\cup\tau\in K \}.
\end{align*}
If $v\in K$ is a vertex, then $\star_K v$ is the subcomplex
consisting of all simplices of $K$ containing $v$, and all their
subsimplices. Note also that $\star_K v$ is the cone over $\link_K
v$.

The following fundamental theorem characterises Cohen--Macaulay
com\-p\-lex\-es combinatorially.

\begin{theo}[Reisner]\label{reisner}
A simplicial complex $K$ is Cohen--Macaulay over $\k$ if and only
if for any simplex $\s\in K$ (including $\s=\varnothing$) and
$i<\dim(\link_K\s)$, it holds that
$\widetilde{H}_i(\link_K\s;\k)=0$.
\end{theo}

Using standard techniques of $PL$ topology the previous theorem
may be reformulated in purely topological terms.

\begin{prop}[Munkres]
  $K^{n-1}$ is Cohen--Macaulay over $\k$ if and only if for
  an 
  arbitrary point $x\in|K|$, it holds that
  \[
    \widetilde{H}_i(|K|;\k)=H_i(|K|,|K|\backslash\,x;\k)=0
    \quad\text{ for }i<n-1.
  \]
\end{prop}

Thus any triangulation of a sphere is a Cohen--Macaulay complex.

\subsection{Resolutions and $\Tor$-algebras}
Let $M$ be a finitely-generated graded
$\k[v_1,\ldots,v_m]$-module. A \emph{free resolution} of $M$ is an
exact sequence
\begin{equation}
\label{resol}
\begin{CD}
  \ldots @>d>> R^{-i} @>d>> \ldots @>d>> R^{-1} @>d>> R^0 @>>> M\to 0,
\end{CD}
\end{equation}
where the $R^{-i}$ are finitely-generated graded free
$\k[v_1,\ldots,v_m]$-modules and the maps $d$ are
degree-preserving. By the Hilbert syzygy theorem, there is a free
resolution of $M$ with $R^{-i}=0$ for $i>m$. A
resolution~(\ref{resol}) determines a bigraded differential
$\k$-module $[R,d]$, where $R=\bigoplus R^{-i,j}$, \
$R^{-i,j}:=(R^{-i})^j$ and $d\colon R^{-i,j}\to R^{-i+1,j}$. The
bigraded cohomology module $H[R,d]$ has $H^{-i,k}[R,d]=0$ for
$i>0$ and $H^{0,k}[R,d]=M^k$. Let $[M,0]$ be the bigraded module
with $M^{-i,k}=0$ for $i>0$, \ $M^{0,k}=M^k$, and zero
differential. Then the resolution~\eqref{resol} determines a
bigraded map $[R,d]\to[M,0]$ inducing an isomorphism in
cohomology.

Let $N$ be another module; then applying the functor
$\otimes_{\k[v_1,\ldots,v_m]}N$ to a resolution $[R,d]$ we get a
homomorphism of differential modules
\[
  [R\otimes_{\k[v_1,\ldots,v_m]}N,d]\to[M\otimes_{\k[v_1,\ldots,v_m]}N,0],
\]
which in general does not induce an isomorphism in cohomology. The
$(-i)$th cohomology module of the cochain complex
\[
\begin{CD}
  \ldots @>>> R^{-i}\otimes_{\k[v_1,\ldots,v_m]}N @>>> \ldots @>>>
  R^0\otimes_{\k[v_1,\ldots,v_m]}N @>>> 0
\end{CD}
\]
is denoted by $\Tor^{-i}_{\k[v_1,\ldots,v_m]}(M,N)$. Thus,
\[
  \Tor^{-i}_{\k[v_1,\ldots,v_m]}(M,N):=
  \frac{\mathop{\mathrm{Ker}}\bigl[d\colon R^{-i}\otimes_{\k[v_1,\ldots,v_m]}
  N\to R^{-i+1}\otimes_{\k[v_1,\ldots,v_m]} N\bigr]}
  {d(R^{-i-1}\otimes_{\k[v_1,\ldots,v_m]}N)}.
\]
Since all the $R^{-i}$ and $N$ are graded modules, we actually
have a \emph{bigraded} $\k$-module
\[
  \Tor_{\k[v_1,\ldots,v_m]}(M,N)=
  \bigoplus_{i,j}\Tor^{-i,j}_{\k[v_1,\ldots,v_m]}(M,N).
\]

The following properties of $\Tor^{-i}_{\k[v_1,\ldots,v_m]}(M,N)$
are well known.
\begin{prop}
\label{torprop} {\rm (a)} the module
$\Tor^{-i}_{\k[v_1,\ldots,v_m]}(M,N)$ does not depend on a choice
of resolution in~{\rm(\ref{resol})};

{\rm (b)} $\Tor^{-i}_{\k[v_1,\ldots,v_m]}(\:\cdot\:,N)$ and
$\Tor^{-i}_{\k[v_1,\ldots,v_m]}(M,\:\cdot\:)$ are covariant
functors;

{\rm (c)} $\Tor^0_{\k[v_1,\ldots,v_m]}(M,N)\cong
M\otimes_{\k[v_1,\ldots,v_m]}N$;

{\rm (d)} $\Tor^{-i}_{\k[v_1,\ldots,v_m]}(M,N)\cong
\Tor^{-i}_{\k[v_1,\ldots,v_m]}(N,M)$.
\end{prop}

Now put $M=\k[K]$ and $N=\k$. Since $\deg v_i=2$, we have
\[
  \Tor_{\k[v_1,\ldots,v_m]}\bigl(\k[K],\k\bigr)=
  \bigoplus_{i,j=0}^m\Tor^{-i,2j}_{\k[v_1,\ldots,v_m]}\bigl(\k[K],\k\bigr)
\]
Define the \emph{bigraded Betti numbers} of $\k[K]$ by
\begin{equation}
\label{bbnfr}
  \b^{-i,2j}\bigl(\k[K]\bigr):=
  \dim_\k\Tor^{-i,2j}_{\k[v_1,\ldots,v_m]}\bigl(\k[K],\k\bigr),\qquad
  0\le i,j\le m.
\end{equation}
We also set
\[
  \b^{-i}(\k[K])=\dim_\k\Tor^{-i}_{\k[v_1,\ldots,v_m]}(\k[K],\k)=
  \sum_j\b^{-i,2j}(\k[K]).
\]

\begin{exam}\label{ressqu}
Let $K$ be the boundary of a square. Then
\[
  \k[K]\cong\k[v_1,\ldots,v_4]/(v_1v_3,v_2v_4).
\]
Let us construct a resolution of $\k[K]$ and calculate the
corresponding bigraded Betti numbers. The module $R^0$ has one
generator 1 (of degree 0), and the map $R^0\to\k[K]$ is the
quotient projection. Its kernel is the ideal $\mathcal I_{K}$,
generated by two monomials $v_1v_3$ and $v_2v_4$. Take $R^{-1}$ to
be a free module on two 4-dimensional generators, denoted $v_{13}$
and $v_{24}$, and define $d\colon R^{-1}\to R^0$ by sending
$v_{13}$ to $v_1v_3$ and $v_{24}$ to $v_2v_4$. Its kernel is
generated by one element $v_2v_4v_{13}-v_1v_3v_{24}$. Hence,
$R^{-2}$ has one generator of degree~8, say~$a$, and the map
$d\colon R^{-2}\to R^{-1}$ is injective and sends $a$ to
$v_2v_4v_{13}-v_1v_3v_{24}$. Thus, we have a resolution
\[
\begin{CD}
  0 @>>> R^{-2} @>>> R^{-1} @>>> R^0 @>>> M @>>> 0
\end{CD}
\]
where $\mathop{\rm rank}R^{0}=\b^{0,0}(\k[K])=1$,\ $\mathop{\rm
rank}R^{-1}=\b^{-1,4}=2$ and $\mathop{\rm
rank}R^{-2}=\b^{-2,8}=1$.
\end{exam}

The Betti numbers $\beta^{-i,2j}(\k[K])$ are important
combinatorial invariants of
the 
simplicial complex~$K$. The following result expresses them in terms
of homology groups of subcomplexes of~$K$.

Given a subset $\omega\subseteq[m]$, we may restrict $K$ to
$\omega$ and consider the \emph{full subcomplex}
$K_\omega=\{\sigma\in K\colon\s\subseteq\omega\}$.

\begin{theo}[Hochster]\label{hoch}
We have
\[
  \beta^{-i,2j}\bigl(\k[K]\bigr)=\sum_{\omega\subseteq[m]\colon|\omega|=j}
  \dim_\k\widetilde{H}^{j-i-1}(K_\omega;\k),
\]
where $\widetilde H^*(\cdot)$ denotes the reduced cohomology
groups and we assume that $\widetilde{H}^{-1}(\varnothing)=\k$.
\end{theo}

Hochster's original proof of this theorem uses rather complicated
combinatorial and commutative algebra techniques. Later in
subsection~\ref{mhoch} we give a topological interpretation of the
numbers $\beta^{-i,2j}(\k[K])$ as the bigraded Betti numbers of a
topological space, and prove a generalisation of Hochster's
theorem.

\begin{exam}[Koszul resolution]
\label{koszul}
Let $M=\k$ with the $\k[v_1,\ldots,v_m]$-module
structure defined via the map $\k[v_1,\ldots,v_m]\to\k$ sending
each $v_i$ to~0. Let $\L[u_1,\ldots,u_m]$ denote the exterior
$\k$-algebra on $m$ generators. The tensor product
$R=\L[u_1,\ldots,u_m]\otimes\k[v_1,\ldots,v_m]$ (here and below we
use $\otimes$ for $\otimes_\k$) may be turned to a differential
bigraded algebra by setting
\begin{gather}
  \bideg u_i=(-1,2),\quad\bideg v_i=(0,2),\notag\\
  \label{diff}
  du_i=v_i,\quad dv_i=0,
\end{gather}
and requiring $d$ to be a derivation of algebras. An explicit
construction of a cochain homotopy shows that $H^{-i}[R,d]=0$ for
$i>0$ and $H^0[R,d]=\k$. Since
$\L[u_1,\ldots,u_m]\otimes\k[v_1,\ldots,v_m]$ is a free
$\k[v_1,\ldots,v_m]$-module, it determines a free resolution
of~$\k$. It is known as the \emph{Koszul resolution} and its
expanded form~\eqref{resol} is as follows:
\begin{multline*}
0\to\L^m[u_1,\ldots,u_m]\otimes\k[v_1,\ldots,v_m]\longrightarrow\cdots\\
\longrightarrow\L^1[u_1,\ldots,u_m]\otimes\k[v_1,\ldots,v_m]\longrightarrow
\k[v_1,\ldots,v_m]\longrightarrow\k\to0
\end{multline*}
where $\L^i[u_1,\ldots,u_m]$ is the subspace of
$\L[u_1,\ldots,u_m]$ spanned by monomials of length~$i$.
\end{exam}

Now let us consider the differential bigraded algebra
$[\L[u_1,\ldots,u_m]\otimes\k[K],d]$ with $d$ defined as
in~{\rm(\ref{diff})}.

\begin{lemm}\label{koscom} There is an isomorphism of bigraded
modules:
\[
  \Tor_{\k[v_1,\ldots,v_m]}(\k[K],\k)\cong
  H\bigl[\L[u_1,\ldots,u_m]\otimes\k[K],d\bigr]
\]
which endows $\Tor_{\k[v_1,\ldots,v_m]}(\k[K],\k)$ with a bigraded
algebra structure in a canonical way.
\end{lemm}
\begin{proof}
Using the Koszul resolution in the definition of $\Tor$, we
calculate
\begin{align*}
    \Tor_{\k[v_1,\ldots,v_m]}(\k[K],\k)
    &\cong\Tor_{\k[v_1,\ldots,v_m]}(\k,\k[K])\\
    &=H\bigl[ \L[u_1,\ldots,u_m]\otimes\k[v_1,\ldots,v_m]
    \otimes_{\k[v_1,\ldots,v_m]}\k[K] \bigr]\\
    &\cong H\bigl[ \L[u_1,\ldots,u_m]\otimes\k[K]\bigr].
\end{align*}
The cohomology in the right hand side is a bigraded algebra,
providing an algebra structure for
$\Tor_{\k[v_1,\ldots,v_m]}(\k[K],\k)$.
\end{proof}

The bigraded algebra $\Tor_{\k[v_1,\ldots,v_m]}(\k[K],\k)$ is
called the \emph{$\Tor$-algebra} of the simplicial complex~$K$.

\begin{lemm}
\label{tamap} A simplicial map $\phi\colon K_1\to K_2$ between two
simplicial complexes on the vertex sets $[m_1]$ and $[m_2]$
respectively induces a homomorphism
\begin{equation}\label{toralgmap}
  \phi_t^*\colon \Tor_{\k[w_1,\ldots,w_{m_2}]}(\k[K_2],\k)\to
  \Tor_{\k[v_1,\ldots,v_{m_1}]}(\k[K_1],\k)
\end{equation}
of the corresponding $\Tor$-algebras.
\end{lemm}
\begin{proof}
This follows directly from Propositions~\ref{frmap}
and~\ref{torprop}~(b).
\end{proof}

\section{Toric spaces}
Moment-angle complexes provide a functor $K\mapsto\zk$ from the
category of simplicial complexes and simplicial maps to the
category of spaces with torus action and equivariant maps. This
functor allows us to use the techniques of equivariant topology in
the study of combinatorics of simplicial complexes and commutative
algebra of their face rings; in a way, it breathes a geometrical
life into Stanley's `combinatorial commutative algebra'. In
particular, the calculation of the cohomology of $\zk$ opens a way
to a topological treatment of homological invariants of face
rings.

The space $\zk$ was introduced for arbitrary finite simplicial
complex $K$ by Davis and Januszkiewicz~\cite{da-ja91} as a
technical tool in their study of (quasi)toric manifolds, a
topological generalisation of smooth algebraic toric varieties.
Later this space turned out to be of great independent interest.
For the subsequent study of $\zk$, its place within `toric
topology', and connections with combinatorial problems we refer
to~\cite{bu-pa02} and its extended Russian version~\cite{bu-pa04}.
Here we review the most important aspects of this study related to
the cohomology of face rings.

\subsection{Moment-angle complexes}
The $m$\emph{-torus} $T^m$ is a product of $m$ circles; we usually
regard it as embedded in $\C^m$ in the standard way:
\[
  T^m=\bigl\{ (z_1,\ldots,z_m)\in\C^m\colon |z_i|=1,\quad i=1,\ldots,m
  \bigr\}.
\]
It is contained in the \emph{unit polydisk}
\[
  (D^2)^m=\{ (z_1,\ldots,z_m)\in\C^m\colon |z_i|\le1,\quad i=1,\ldots,m
  \}.
\]
For
an 
arbitrary subset $\omega\subseteq V$, define
\[
  B_\omega:=\{(z_1,\ldots,z_m)\in(D^2)^m\colon |z_i|=1\text{ for
  }i\notin\omega\}.
\]
The subspace $B_\omega$ is homeomorphic to $(D^2)^{|\omega|}\times
T^{m-|\omega|}$.

Given a simplicial complex $K$ on $[m]=\{1,\ldots,m\}$, we define
the \emph{moment-angle complex} $\zk$ by
\begin{equation}\label{zkuni}
  \zk:=\bigcup_{\sigma\in K}B_\sigma\subseteq(D^2)^m.
\end{equation}

The torus $T^m$ acts on $(D^2)^m$ coordinatewise and each subspace
$B_\omega$ is invariant under this action. Therefore, the space
$\zk$ inherits a torus action. The quotient $(D^2)^m/T^m$ can be
identified with the \emph{unit $m$-cube}:
\[
  I^m:=\bigl\{ (y_1,\ldots,y_m)\in\R^m\colon 0\le y_i\le1,\quad i=1,\ldots,m
  \bigr\}.
\]
The quotient $B_\omega/T^m$ is then the following
$|\omega|$-dimensional face of $I^m$:
\[
  C_\omega:=\bigl\{ (y_1,\ldots,y_m)\in I^m\colon y_i=1\text{ if
  }i\notin\omega\bigr\}.
\]
Thus the whole quotient $\zk/T^m$ is identified with a certain
cubical subcomplex in $I^m$, which we denote by $\cc(K)$.

\begin{lemm}
The cubical complex $\cc(K)$ is $PL$-homeomorphic to $\cone K$.
\end{lemm}
\begin{proof}
Let $K'$ denote the barycentric subdivision of $K$ (the vertices
of $K'$ correspond to non-empty simplices $\sigma$ of $K$). We
define a $PL$ embedding $i_c\colon \cone K'\hookrightarrow I^m$ by
mapping each vertex $\sigma$ to the vertex $(\e_1,\dots,\e_m)\in
I^m$ where $\e_i=0$ if $i\in\sigma$ and $\e_i=1$ otherwise, the
cone vertex to $(1,\dots,1)\in I^m$, and then extending linearly
on the simplices of $\cone K'$. The barycentric subdivision of a
face $\sigma\in K$ is a subcomplex in $K'$, which we denote
$K'|_\sigma$. Under the map $i_c$ the subcomplex $\cone
K'|_\sigma$ maps onto the face $C_\sigma\subset I^m$. Thus the
whole complex $\cone K'$ maps homeomorphically onto $\cc(K)$,
which concludes the proof.
\end{proof}

It follows that the moment-angle complex $\zk$ can be defined by
the pullback diagram
\[
\begin{CD}
  \zk @>>> (D^2)^m\\
  @VVV @VV\rho V\\
  \cone K' @>i_c>> I^m
\end{CD}
\]
where $\rho$ is the projection onto the orbit space.

\begin{exam}
The embedding $i_c$ for two simple cases when $K$ is a three point
complex and the boundary of a triangle is shown on
Figure~\ref{cubem}. If $K=\D^{m-1}$ is the whole simplex on $m$
vertices, then $\cc(K)$ is the whole cube $I^m$, and the above
constructed $PL$-homeomorphism between $\cone(\D^{m-1})'$ and
$I^m$ defines the \emph{standard triangulation} of~$I^m$.
\begin{figure}
  \begin{picture}(120,45)
  \put(15,5){\line(1,0){25}}
  \put(15,5){\line(0,1){25}}
  \put(40,5){\line(1,1){10}}
  \put(50,15){\line(0,1){25}}
  \put(50,40){\line(-1,0){25}}
  \put(25,40){\line(-1,-1){10}}
  \put(40,5){\circle*{2}}
  \put(40,30){\circle*{2}}
  \put(15,30){\circle*{2}}
  \put(50,40){\circle*{2}}
  \multiput(40,29.3)(0,0.1){16}{\line(1,1){10}}
  \multiput(25,15)(5,0){5}{\line(1,0){3}}
  \multiput(25,15)(0,5){5}{\line(0,1){3}}
  \multiput(25,15)(-5,-5){2}{\line(-1,-1){4}}
  \put(26,16){$0$}
  \put(22,-2){$K=\text{3 points}$}
  \put(75,5){\line(1,0){25}}
  \put(75,5){\line(0,1){25}}
  \put(100,5){\line(1,1){10}}
  \put(110,15){\line(0,1){25}}
  \put(110,40){\line(-1,0){25}}
  \put(85,40){\line(-1,-1){10}}
  \put(100,5){\circle*{2}}
  \put(110,15){\circle*{2}}
  \put(100,30){\circle*{2}}
  \put(75,30){\circle*{2}}
  \put(85,40){\circle*{2}}
  \put(75,5){\circle*{2}}
  \put(110,40){\circle*{2}}
  \multiput(100,29.3)(0,0.1){16}{\line(1,1){10}}
  \multiput(75,29.3)(0,0.1){16}{\line(1,1){10}}
  \multiput(100,4.3)(0,0.1){16}{\line(1,1){10}}
  \multiput(85,15)(5,0){5}{\line(1,0){3}}
  \multiput(85,15)(0,5){5}{\line(0,1){3}}
  \multiput(85,15)(-5,-5){2}{\line(-1,-1){4}}
  \multiput(78,5)(4,0){6}{\line(0,1){25}}
  \multiput(78,30)(4,0){6}{\line(1,1){10}}
  \multiput(100,8)(0,4){6}{\line(1,1){10}}
  \put(86,16){$0$}
  \put(82,-2){$K=\partial\varDelta^2$}
  \linethickness{1mm}
  \put(15,30){\line(1,0){25}}
  \put(40,5){\line(0,1){25}}
  \put(75,30){\line(1,0){25}}
  \put(75,5){\line(1,0){25}}
  \put(85,40){\line(1,0){25}}
  \put(100,5){\line(0,1){25}}
  \put(75,5){\line(0,1){25}}
  \put(110,15){\line(0,1){25}}
  \end{picture}
\caption{Embedding $i_c\colon\cone K'\hookrightarrow I^m$.}
\label{cubem}
\end{figure}
\end{exam}

The next lemma shows that the space $\zk$ is particularly nice for
certain geometrically important classes of triangulations.

\begin{lemm}\label{zkman}
Suppose that $K$ is a triangulation of an $(n-1)$-dimensional
sphere. Then $\zk$ is a closed $(m+n)$-dimensional manifold.

{\sloppy In general, if $K$ is a triangulated manifold then
$\zk\setminus\rho^{-1}(1,\ldots,1)$ is a noncompact manifold,
where $(1,\ldots,1)\in I^m$ is the cone vertex and
$\rho^{-1}(1,\ldots,1)\cong T^m$.

}
\end{lemm}
\begin{proof}
We only prove the first statement here; the proof of the second is
similar and can be found in~\cite{bu-pa04}. Each vertex $v_i$ of
$K$ corresponds to a vertex of the barycentric subdivision $K'$,
which we continue to denote~$v_i$. Let $\star_{K'}v_i$ be the star
of $v_i$ in $K'$, that is, the subcomplex consisting of all
simplices of $K'$ containing $v_i$, and all their subsimplices.
The space $\cone K'$ has a canonical \emph{face structure} whose
facets (co\-di\-men\-si\-on-\-one faces) are
\begin{equation}\label{facet}
  F_i:=\star_{K'}v_i, \quad i=1,\ldots,m,
\end{equation}
and whose $i$-faces are non-empty intersections of $i$-tuples of
facets. In particular, the vertices (0-faces) in this face
structure are the barycentres of $(n-1)$-dimensional simplices
of~$K$.

For every such barycentre $b$ we denote by $U_b$ the subset of
$\cone K'$ obtained by removing all faces not containing $b$.
Since $K$ is a triangulation of a sphere, $\cone K'$ is an
$n$-ball, hence each $U_b$ is homeomorphic to an open subset in
$I^n$ via a homeomorphism preserving the dimension of faces. Since
each point of $\cone K'$ is contained in some $U_b$, this displays
$\cone K'$ as a \emph{manifold with corners}. Having identified
$\cone K'$ with $\cc(K)$ and further $\cc(K)$ with $\zk/T^m$, we
see that every point in $\zk$ lies in a neighbourhood homeomorphic
to an open subset in $(D^2)^n\times T^{m-n}$ and thus in
$\R^{m+n}$.
\end{proof}

A particularly important class of examples of sphere
triangulations arise from boundary triangulations of convex
polytopes. Suppose $P$ is a \emph{simple} $n$-dimensional convex
polytope, i.e. one where every vertex is contained in exactly $n$
facets. Then the dual (or \emph{polar}) polytope is
\emph{simplicial}, and we denote its boundary complex by $K_P$.
$K_P$ is then a triangulation of an $(n-1)$-sphere. The faces of
$\cone K'_P$ introduced in the previous proof coincide with those
of~$P$.

\begin{exam}
Let $K=\partial\Delta^{m-1}$. Then $\zk=\partial((D^2)^m)\cong
S^{2m-1}$. In particular, for $m=2$ from~\eqref{zkuni} we get the
familiar decomposition
\[
  S^3=D^2\times S^1\cup S^1\times D^2\subset D^2\times D^2
\]
of a 3-sphere into a union of two solid tori.
\end{exam}

Using faces~\eqref{facet} we can identify the isotropy subgroups
of the $T^m$-action on $\zk$. Namely, the isotropy subgroup of a
point $x$ in the orbit space $\cone K'$ is the coordinate subtorus
\[
  T(x)=\{ (z_1,\ldots,z_m)\in T^m\colon z_i=1\text{ if }x\notin F_i
  \}.
\]
In particular, the action is free over the interior (that is, near
the cone point) of $\cone K'$.

It follows that the moment-angle complex can be identified with
the quotient
\[
  \zk=\bigl(T^m\times|\cone K'|\bigr)/{\sim},
\]
where $(t_1,x)\sim(t_2,y)$ if and only if $x=y$ and
$t_1t_2^{-1}\in T(x)$. In the case when $K$ is the dual
triangulation of a simple polytope $P^n$ we may write $(T^m\times
P^n)/{\sim}$ instead. The latter $T^m$-manifold is the one
introduced by Davis and Januszkiewicz~\cite{da-ja91}, which
thereby coincides with our moment-angle complex.

\subsection{Homotopy fibre construction}
The classifying space for the circle $S^1$ can be identified with
the infinite-dimensional projective space $\C P^\infty$. The
classifying space $BT^m$ of the $m$-torus is a product of $m$
copies of $\C P^\infty$. The cohomology of $BT^m$ is the
polynomial ring $\Z[v_1,\ldots,v_m]$, \ $\deg v_i=2$ (the
cohomology is taken with integer coefficients, unless another
coefficient ring is explicitly specified). The total space $ET^m$
of the universal principal $T^m$-bundle over $BT^m$ can be
identified with the product of $m$ infinite-dimensional spheres.

In~\cite{da-ja91} Davis and Januszkiewicz considered the
\emph{homotopy quotient} of $\zk$ by the $T^m$-action (also known
as the \emph{Borel construction}). We refer to it as the
\emph{Davis--Januszkiewicz space}:
\[
  \djs(K):=ET^m\times_{T^m}\zk=ET^m\times\zk/\!\sim,
\]
where $(e,z)\sim(et^{-1},tz)$. There is a a fibration $p\colon
\djs(K)\to BT^m$ with fibre~$\zk$.
The cohomology 
of the Borel construction of a $T^m$-space $X$ is called the
\emph{equivariant cohomology} and denoted by $H^*_{T^m}(X)$.

A theorem of~\cite{da-ja91} states that the cohomology ring of
$\djs(K)$ (or the equivariant cohomology of $\zk$) is isomorphic
to $\Z[K]$. This result can be clarified by an alternative
construction of $\djs(K)$~\cite{bu-pa02}, which we review below.

The space $BT^m$ has the canonical cell decomposition in which
each factor $\C P^\infty$ has one cell in every even dimension.
Given a subset $\omega\subseteq[m]$, define the subproduct
\[
  BT^\omega:=\{(x_1,\dots,x_m)\in BT^m\colon x_i=*\text{ if
  }i\notin\omega\}
\]
where $*$ is the basepoint (zero-cell) of~$\C P^{\infty}$. Now for
a simplicial complex $K$ on $[m]$ define the following cellular
subcomplex:
\begin{equation}\label{djuni}
  BT^K:=\bigcup_{\s\in K}BT^\s \subseteq BT^m.
\end{equation}

\begin{prop}
\label{homsrs} The cohomology of $BT^K$ is isomorphic to the
Stanley--Reisner ring~$\Z[K]$. Moreover, the inclusion of cellular
complexes $i\colon BT^K\hookrightarrow BT^m$ induces the quotient
epimorphism
\[
  i^*\colon \Z[v_1,\ldots,v_m]\to\Z[K]=\Z[v_1,\ldots,v_m]/\mathcal I_K
\]
in the cohomology.
\end{prop}
\begin{proof}
Let $B^{2k}_i$ denote the $2k$-dimensional cell in the $i$th
factor of $BT^m$, and $C^*(BT^m)$ the cellular cochain module. A
monomial $v_{i_1}^{k_1}\dots v_{i_p}^{k_p}$ represents the
cellular cochain $(B^{2k_1}_{i_1}\dots B^{2k_p}_{i_p})^*$ in
$C^*(BT^m)$. Under the cochain homomorphism induced by the
inclusion $BT^K\subset BT^m$ the cochain $(B^{2k_1}_{i_1}\dots
B^{2k_p}_{i_p})^*$ maps identically if $\{i_1,\dots,i_p\}\in K$
and to zero otherwise, whence the statement follows.
\end{proof}

\begin{theo}
There is a deformation retraction $\djs(K)\to BT^K$ such that the
diagram
\[
\begin{CD}
  \djs(K) @>p>> BT^m\\
  @VVV @|\\
  BT^K @>i>> BT^m
\end{CD}
\]
is commutative.
\end{theo}
\begin{proof}
We have $\zk=\bigcup_{\s\in K}B_\s$, and each $B_\s$ is
$T^m$-invariant. Hence, there is the corresponding decomposition
of the Borel construction:
\[
  \djs(K)=ET^m\times_{T^m}\zk=\bigcup_{\s\in K}ET^m\times_{T^m}B_\s.
\]
Suppose $|\s|=s$. Then $B_\s\cong(D^2)^s\times T^{m-s}$, so we
have
\[
  ET^m\times_{T^m}B_\s\cong(ET^s\times_{T^s}(D^2)^s)\times
  ET^{m-s}.
\]
The space $ET^s\times_{T^s}(D^2)^s$ is the total space of a
$(D^2)^s$-bundle over $BT^s$, and $ET^{m-s}$ is contractible. It
follows that there is a deformation retraction
$ET^m\times_{T^m}B_\s\to BT^\s$. These homotopy equivalences
corresponding to different simplices fit together to yield
the 
required homotopy equivalence between $p\colon\djs(K)\to BT^m$ and
$i\colon BT^K\hookrightarrow BT^m$.
\end{proof}

\begin{coro}\label{hofib}
The space $\zk$ is the homotopy fibre of the cellular inclusion
$i\colon BT^K\hookrightarrow BT^m$. Hence~\cite{da-ja91} there are
ring isomorphisms
\[
  H^*(\djs(K))=H^*_{T^m}(\zk)\cong\Z[K].
\]
\end{coro}

In view of the last two statements we shall also use the notation
$\djs(K)$ for $BT^K$, and refer to the whole class of spaces
homotopy equivalent to $\djs(K)$ as the \emph{Davis--Januszkiewicz
homotopy type}.

An important question arises: to what extent does the isomorphism
of the cohomology ring of a space $X$ with the face ring $\Z[K]$
determine the homotopy type of $X$? In other words, for given $K$,
does there exist a `fake' Davis--Januszkiewicz space, whose
cohomology is isomorphic to $\Z[K]$, but which is not homotopy
equivalent to $\djs(K)$? This question is addressed
in~\cite{no-ra05}. It is shown there~\cite[Prop.~5.11]{no-ra05}
that if $\Q[K]$ is a \emph{complete intersection ring} and $X$ is
a nilpotent cell complex of finite type whose rational cohomology
is isomorphic to $\Q[K]$, then $X$ is rationally homotopy
equivalent to $\djs(K)$. Using the formality of~$\djs(K)$, this
can be rephrased by saying that the complete intersection face
rings are \emph{intrinsically formal} in the sense of Sullivan.

Note that the class of simplicial complexes $K$ for which the face
ring $\Q[K]$ is a complete intersection has a transparent
geometrical interpretation: such $K$ is a join of simplices and
boundaries of simplices.

\subsection{Coordinate subspace arrangements}
Yet another interpretation of the moment-angle complex $\zk$ comes
from its identification up to homotopy with the complement of the
complex coordinate subspace arrangement corresponding to~$K$. This
leads to an application of toric topology in the theory of
arrangements, and allows us to describe and effectively calculate
the cohomology rings of coordinate subspace arrangement
complements and in certain cases identify their homotopy types.

A \emph{coordinate subspace} in $\C^m$ can be written as
\begin{equation}\label{cossp}
  L_\omega=\{(z_1,\ldots,z_m)\in\C^m\colon
  z_{i_1}=\cdots=z_{i_k}=0\}
\end{equation}
for some subset $\omega=\{i_1,\ldots,i_k\}\subseteq[m]$. Given a
simplicial complex $K$, we may define the corresponding
\emph{coordinate subspace arrangement}
$\{L_\omega\colon\omega\notin K\}$ and its \emph{complement}
\[
  U(K)=\C^m\setminus\bigcup_{\omega\notin K}L_\omega.
\]
Note that if $K'\subset K$ is a subcomplex, then $U(K')\subset
U(K)$. It is easy to see~\cite[Prop.~8.6]{bu-pa02} that the
assignment $K\mapsto U(K)$ defines a one-to-one order preserving
correspondence between the set of simplicial complexes on $[m]$
and the set of coordinate subspace arrangement complements
in~$\C^m$.

The subset $U(K)\subset\C^m$ is invariant with respect to the
coordinatewise $T^m$-action. It follows from~\eqref{zkuni} that
$\zk\subset U(K)$.

\begin{prop}\label{ukzk}
There is a $T^m$-equivariant deformation retraction
\[
  U(K)\stackrel\simeq\longrightarrow\zk.
\]
\end{prop}
\begin{proof}
In analogy with~\eqref{djuni}, we may write
\begin{equation}\label{ukuni}
  U(K)=\bigcup_{\s\in K}U_\s,
\end{equation}
where
\[
  U_\sigma:=\{(z_1,\ldots,z_m)\in\C^m\colon z_i\ne0\text{ for
  }i\notin\sigma\}.
\]
Then there are obvious homotopy equivalences (deformation
retractions)
\[
  \C^\s\times(\C\setminus0)^{[m]\setminus\s}\cong U_\s
  \stackrel\simeq\longrightarrow B_\s\cong
  (D^2)^\s\times(S^1)^{[m]\setminus\s}.
\]
These patch together to get the required map $U(K)\to\zk$.
\end{proof}

\begin{exam}\label{csaex}
\begin{enumerate}
\item Let $K=\partial\varDelta^{m-1}$. Then $U(K)=\C^m\setminus0$
(recall that $\zk\cong S^{2m-1}$ in this case).

\item Let $K=\{v_1,\ldots,v_m\}$ ($m$ points). Then
\[
  U(K)=\C^m\setminus\bigcup_{1\le i<j\le m}\{z_i=z_j=0\},
\]
the complement to the set of all codimension 2 coordinate planes.

\item More generally, if $K$ is the $i$-skeleton of
$\varDelta^{m-1}$, then $U(K)$ is the complement to the set of all
coordinate planes of codimension $(i+2)$.
\end{enumerate}
\end{exam}

The reader may have noticed a similar pattern in several
constructions of toric spaces appeared above;
compare~\eqref{zkuni},~\eqref{djuni} and~\eqref{ukuni}. The
following general framework was suggested to the author by Neil
Strickland in a private communication.

\begin{cons}[$K$-power]\label{nsc}
Let $X$ be a space and $W\subset X$ a subspace. For a simplicial
complex $K$ on $[m]$ and $\s\in K$, we set
\begin{align*}
  (X,W)^\s:=\bigl\{(x_1,\ldots,x_m)\in
  X^m\colon x_j\in W\text{ for }j\notin\s\bigl\}\\
\intertext{and}
  (X,W)^K:=\bigcup_{\s\in K}(X,W)^\s=
  \bigcup_{\s\in K}\Bigl(\prod_{i\in\s}X\times\prod_{i\notin\s}W \Bigl).
\end{align*}
We refer to the space $(X,W)^K\subseteq X^m$ as the
\emph{$K$-power} of $(X,W)$. If $X$ is a pointed space and $W=pt$
is the basepoint, then we shall use the abbreviated notation
$X^K:=(X,pt)^K$. Examples considered above include
$\zk=(D^2,S^1)^K$, $\cc(K)=(I^1,S^0)^K$, $\djs(K)=(\C P^\infty)^K$
and $U(K)=(\C,\C^*)^K$.
\end{cons}

Homotopy theorists would recognise the $K$-power as an example of
the \emph{colimit} of a diagram of topological spaces over the
\emph{face category} of $K$ (objects are simplices and morphisms
are inclusions). The diagram assigns the space $(X,W)^\sigma$ to a
simplex $\sigma$; its colimit is $(X,W)^K$. These observations are
further developed and used to construct models of loop spaces of
toric spaces as well as for homotopy and homology calculations
in~\cite{p-r-v04} and~\cite{pa-ra??}.

\subsection{Toric varieties, quasitoric manifolds, and torus manifolds}
Several important classes of manifolds with torus action emerge as
the quotients of moment-angle complexes by appropriate freely
acting subtori.

First we give the following characterisation of lsops in the face
ring. Let $K^{n-1}$ be a simplicial complex and $t_1,\dots,t_n$ a
sequence of degree-two elements in $\k[K]$. We may write
\begin{equation}\label{lsop}
  t_i=\l_{i1}v_1+\dots+\l_{im}v_m,\quad i=1,\dots,n.
\end{equation}
For
an 
arbitrary simplex $\s\in K$, we have $K_\s=\D^{|\s|-1}$ and
$\k[K_\s]$ is the polynomial ring $\k[v_i\colon i\in\s]$ on $|\s|$
generators. The inclusion $K_\s\subset K$ induces the
\emph{restriction homomorphism} $r_\s$ from $\k[K]$ to the
polynomial ring, mapping $v_i$ identically if $i\in\s$ and to zero
otherwise.

\begin{lemm}\label{restr}
A degree-two sequence $t_1,\ldots,t_n$ is an lsop in $\k[K^{n-1}]$
if and only if for every $\s\in K$ the elements
$r_\s(t_1),\ldots,r_\s(t_n)$ generate the positive ideal
$\k[v_i\colon i\in\s]_+$.
\end{lemm}
\begin{proof}
Suppose \eqref{lsop} is an lsop. For simplicity we denote its
image under any restriction homomorphism by the same letters. Then
the restriction induces an epimorphism of the quotient rings:
\[
  \k[K]/(t_1,\dots,t_n)\to\k[v_i\colon i\in\s]/(t_1,\dots,t_n).
\]
Since \eqref{lsop} is an lsop, $\k[K]/(t_1,\dots,t_n)$ is a
finitely generated $\k$-module. Hence, so is $\k[v_i\colon
i\in\s]/(t_1,\dots,t_n)$. But the latter can be finitely generated
only if $t_1,\dots,t_n$ generates $\k[v_i\colon i\in\s]_+$.

The ``if'' part may be proved by considering the sum of
restrictions:
\[
  \k[K]\to\bigoplus_{\s\in K}\k[v_i\colon i\in\s],
\]
which turns out to be a monomorphism. See
\cite[Th.~5.1.16]{br-he98} for details.
\end{proof}

Obviously, it is enough to consider only restrictions to the
maximal simplices in the previous lemma.

Suppose now that $K$ is Cohen--Macaulay (e.g. $K$ is a sphere
triangulation). Then every lsop is a regular sequence (however,
for $\k=\Z$ or a field of finite characteristic an lsop may fail
to exist).

Now we restrict to the case $\k=\Z$ and organise the coefficients
in \eqref{lsop} into an $n\times m$-matrix $\L=(\l_{ij})$. For
an 
arbitrary maximal simplex $\s\in K$ denote by $\L_\s$ the square
submatrix formed by the elements $\l_{ij}$ with $j\in\s$. The
matrix $\Lambda$ defines a linear map $\Z^m\to\Z^n$ and a
homomorphism $T^m\to T^n$. We denote both by~$\l$ and denote the
kernel of the latter map by~$T_\L$.

\begin{theo}\label{zkquo}
The following conditions are equivalent:
\begin{itemize}\label{quot}
\item[(a)] the sequence \eqref{lsop} is an lsop in $\Z[K^{n-1}]$;
\item[(b)] $\det\L_\s=\pm1$ for every maximal simplex $\s\in K$;
\item[(c)] $T_\L\cong T^{m-n}$ and $T_\L$ acts freely on $\zk$.
\end{itemize}
\end{theo}
\begin{proof}
The equivalence of (a) and (b) is a reformulation of
Lemma~\ref{restr}. Let us prove the equivalence of (b) and (c).
Every isotropy subgroup of the $T^m$-action on $\zk$ has the form
\[
  T^\s=\bigl\{ (z_1,\dots,z_m)\in T^m\colon z_i=1\text{ if }i\notin\s \bigr\}
\]
for some simplex $\s\in K$. Now, (b) is equivalent to the
condition $T_\L\cap T^\s=\{e\}$ for arbitrary maximal $\s$, whence
the statement follows.
\end{proof}

We denote the quotient $\zk/T_\L$ by $M^{2n}_K(\L)$, and
abbreviate it to $M^{2n}_K$ or to $M^{2n}$ when the context
allows. If $K$ is a triangulated sphere, then $\zk$ is a manifold,
hence, so is $M^{2n}_K$. The $n$-torus $T^n=T^m/T_\L$ acts on
$M^{2n}_K$. This construction produces two important classes of
$T^n$-manifolds as particular examples.

Let $K=K_P$ be a polytopal triangulation, dual to the boundary
complex of a simple polytope~$P$. Then the map $\l$ determined by
the matrix $\L$ may be regarded as an assignment of an integer
vector to every facet of~$P$. The map $\l$ coming from a matrix
satisfying the condition of Theorem~\ref{quot}(b) was called a
\emph{characteristic map} by Davis and
Januszkiewicz~\cite{da-ja91}. We refer to the corresponding
quotient $M^{2n}_P(\L)=\mathcal Z_{K_P}/T_\L$ as a
\emph{quasitoric manifold} (a toric manifold in the terminology of
Davis--Januszkiewicz).

Let us assume further that $P$ is realised in $\R^n$ with integer
coordinates of vertices, so we can write
\begin{equation}\label{ptope}
  P^n=\bigl\{\mb x\in\R^n\colon\langle\mb l_i,\mb x\rangle\ge-a_i,
  \; i=1,\dots,m\bigr\},
\end{equation}
where $\mb l_i$ are inward pointing normals to the facets of $P^n$
(we may further assume these vectors to be primitive), and
$a_i\in\Q$. Let $\Lambda$ be the matrix formed by the column
vectors $\mb l_i$, \ $i=1,\dots,m$. Then $\mathcal Z_{K_P}/T_\L$
can be identified with the \emph{projective toric
variety}~\cite{dani78,fult93} determined by the polytope~$P$. The
condition of Theorem~\ref{quot}(b) is equivalent to the
requirement that the toric variety is non-singular. Thereby a
non-singular projective toric variety is a quasitoric manifold
(but there are many quasitoric manifolds which are not toric
varieties).

We also note that smooth projective toric varieties provide
examples of \emph{symplectic} $2n$-dimensional manifolds with
\emph{Hamiltonian} $T^n$-action. These symplectic manifolds can be
obtained via the process of \emph{symplectic reduction} from the
standard Hamiltonian $T^m$-action on $\C^m$. A choice of an
$(m-n)$-dimensional toric subgroup provides a \emph{moment map}
$\mu\colon\C^m\to\R^{m-n}$, and the corresponding moment-angle
complex $\mathcal Z_{K_P}$ can be identified with the level
surface $\mu^{-1}(a)$ of the moment map for any of its regular
values~$a$. The details of this construction can be found
in~\cite[p.~130]{bu-pa02}.

Finally, we mention that if $K$ is an arbitrary (not necessarily
polytopal) triangulation of sphere, then the manifold
$M^{2n}_K(\L)$ is a \emph{torus manifold} in the sense of
Hattori--Masuda~\cite{ha-ma03}. The corresponding multi-fan has
$K$ as the underlying simplicial complex. This particular class of
torus manifolds has many interesting properties.

\section{Cohomology of moment-angle complexes}
The main result of this section (Theorem~\ref{zkcoh}) identifies
the \emph{integral} cohomology algebra of the moment-angle complex
$\mathcal Z_K$ with the $\operatorname{Tor}$-algebra of the face
ring of
the 
simplicial complex~$K$. Over the rationals this result was proved
in~\cite{bu-pa99} by studying the Eilenberg--Moore spectral
sequence of the fibration $\zk\to\djs(K)\to BT^m$; a more detailed
account of applications of the Eilenberg--Moore spectral sequence
to toric topology can be found in~\cite{bu-pa02}. The new proof,
which works with integer coefficients as well, relies upon a
construction of a special cellular decomposition of $\zk$ and
subsequent analysis of the corresponding cellular cochains.

One of the key ingredients here is a specific cellular
approximation of the diagonal map $\Delta\colon\mathcal
Z_K\to\mathcal Z_K\times\mathcal Z_K$. Cellular cochains do not
admit a functorial associative multiplication because a proper
cellular diagonal approximation does not exist in general. The
construction of moment-angle complexes is given by a functor from
the category of simplicial complexes to the category of spaces
with a torus action. We show that in this special case the
cellular approximation of the diagonal is functorial with respect
to those maps of moment-angle complexes which are induced by
simplicial maps. The corresponding cellular cochain algebra is
isomorphic to a quotient of the Koszul complex for $\k[K]$ by an
acyclic ideal, and its cohomology is isomorphic to the
$\Tor$-algebra. The proofs have been sketched in~\cite{b-b-p04};
here we follow the more detailed exposition of~\cite{bu-pa04}.
Another proof of Theorem~\ref{zkcoh} follows from a recent
independent work of M.~Franz~\cite[Th.~1.2]{fran06}.

\subsection{Cell decomposition}
The polydisc~$(D^2)^m$ has a cell decomposition in which each
$D^2$ is subdivided into cells 1, $T$ and $D$ of dimensions $0$,
$1$ and~$2$ respectively, see Figure~\ref{cellfig}.
\begin{figure}[h]
  \begin{picture}(120,35)
  \put(60,20){\oval(30,30)}
  \put(75,20){\circle*{1.5}}
  \put(76,21){1}
  \put(42,21){$T$}
  \put(52,25){$D$}
  \end{picture}
  \caption{}
  \label{cellfig}
\end{figure}
Each cell of this complex is a product of cells of 3 different
types and we encode it by a word $\mathcal T\in\{D,T,1\}^m$ in a
three-letter alphabet. Assign to each pair of subsets
$\sigma,\omega\subseteq[m]$, $\sigma\cap\omega=\varnothing$, the
word $\mathcal T(\sigma,\omega)$ which has the letter $D$ on the
positions indexed by $\sigma$ and letter $T$ on the positions with
indices from~$\omega$.

\begin{lemm}
$\zk$ is a cellular subcomplex of $(D^2)^m$. A cell $\mathcal
T(\s,\omega)\subset(D^2)^m$ belongs to $\zk$ if and only if $\s\in
K$.
\end{lemm}
\begin{proof}
We have $\zk=\cup_{\s\in K}B_\s$ and each $B_\s$ is the closure of
the cell $\mathcal T(\s,[m]\setminus\s)$.
\end{proof}

Therefore, we can consider the cellular cochain complex
$C^*(\mathcal Z_K)$, which has an additive basis consisting of the
cochains $\mathcal T(\sigma,\omega)^*$. It has a natural bigrading
defined by
\[
  \bideg\mathcal T(\s,\omega)^*=(-|\omega|,2|\s|+2|\omega|),
\]
so $\bideg D=(0,2)$, $\bideg T=(-1,2)$ and $\bideg 1=(0,0)$.
Moreover, since the cellular differential does not change the
second grading, $C^*(\zk)$ splits into the sum of its components
having fixed second degree:
\[
  C^*(\zk)=\bigoplus_{j=1}^m C^{*,2j}(\zk).
\]
The cohomology 
of $\zk$ thereby acquires an additional grading, and we may define the
\emph{bigraded Betti numbers} $b^{-i,2j}(\zk)$ by
\[
  b^{-i,2j}(\zk):=\mathop{\mathrm{rank}} H^{-i,2j}(\zk), \quad i,j=1,\dots,m.
\]
For the ordinary Betti numbers we have
$b^k(\zk)=\sum_{2j-i=k}b^{-i,2j}(\zk)$.

\begin{lemm}
\label{mamap} Let $\phi\colon K_1\to K_2$ be a simplicial map
between simplicial complexes on the sets $[m_1]$ and $[m_2]$
respectively. Then there is an equivariant cellular map
$\phi_{\mathcal Z}\colon\mathcal Z_{K_1}\to\mathcal Z_{K_2}$
covering the induced map $|\cone K'_1|\to|\cone K'_2|$.
\end{lemm}
\begin{proof}
Define a map of polydisks
\[
  \phi_D\colon (D^2)^{m_1}\to(D^2)^{m_2}, \quad
  (z_1,\ldots,z_{m_1})\mapsto(w_1,\ldots,w_{m_2}),
\]
where
\[
  w_j:=\prod_{i\in\phi^{-1}(j)}z_i,\qquad j=1,\ldots,m_2
\]
(we set $w_j=1$ if $\phi^{-1}(j)=\varnothing$). Assume $\tau\in
K_1$. In the notation of~\eqref{zkuni}, we have
$\phi_D(B_\tau)\subseteq B_{\phi(\tau)}$. Since $\phi$ is a
simplicial map, $\phi(\tau)\in K_2$ and $B_{\phi(\tau)}\subset
\mathcal Z_{K_2}$, so the restriction of $\phi_D$ to $\mathcal
Z_{K_1}$ is the required map.
\end{proof}

\begin{coro}
The correspondence $K\mapsto\zk$ gives rise to a functor from the
category of simplicial complexes and simplicial maps to the
category of spaces with torus actions and equivariant maps. It
induces a natural transformation between the simplicial cochain
functor of $K$ and the cellular cochain functor of $\zk$.
\end{coro}

We also note that the maps respect the bigrading, so the bigraded
Betti numbers are also functorial.

\subsection{Koszul algebras}
Our algebraic model for the cellular cochains of $\zk$ is obtained
by taking the quotient of the Koszul algebra
$[\L[u_1,\ldots,u_m]\otimes\k[K],d]$ from Lemma~\ref{koscom} by a
certain acyclic ideal. Namely, we introduce a factor algebra
\[
  R^*(K):=\L[u_1,\ldots,u_m]\otimes\Z[K]\bigr/(v_i^2=u_iv_i=0,\;
  i=1,\ldots,m),
\]
where the differential and bigrading are as in~\eqref{diff}. Let
\[
  \varrho\colon\L[u_1,\ldots,u_m]\otimes\Z[K]\to R^*(K)
\]
be the quotient projection. The algebra $R^*(K)$ has a finite
additive basis consisting of the monomials of the form $u_\omega
v_\s$ where $\omega\subseteq[m]$, $\s\in K$ and
$\omega\cap\s=\varnothing$ (remember that we are using the
notation $u_\omega=u_{i_1}\ldots u_{i_k}$ for
$\omega=\{i_1,\ldots,i_k\}$). Therefore, we have an additive
inclusion (a monomorphism of bigraded differential modules)
\[
  \iota\colon R^*(K)\to\L[u_1,\ldots,u_m]\otimes\Z[K]
\]
which satisfies $\varrho\cdot\iota=\id$.

The following statement shows that the finite-dimensional quotient
$R^*(K)$ has the same cohomology as the Koszul algebra.

\begin{lemm}\label{iscoh}
The quotient map
$\varrho\colon\L[u_1,\ldots,u_m]\otimes\Z[K]\to R^*(K)$ induces an
isomorphism in cohomology.
\end{lemm}
\begin{proof}
The argument is similar to that used in the proof of the
acyclicity of the Koszul resolution. We construct a cochain
homotopy between the maps $\id$ and $\iota\cdot\varrho$ from
$\L[u_1,\ldots,u_m]\otimes\Z[K]$ to itself, that is, a map $s$
satisfying
\begin{equation}\label{chaineq1}
  ds+sd=\id-\iota\cdot\varrho.
\end{equation}

First assume that $K=\D^{m-1}$. We denote the corresponding
bigraded algebra $\L[u_1,\ldots,u_m]\otimes\Z[K]$ by
\begin{equation}\label{emalg}
  E=E_m:=\L[u_1,\ldots,u_m]\otimes\Z[v_1,\ldots,v_m],
\end{equation}
while $R^*(K)$ is isomorphic to
\begin{equation}\label{rmalg}
  \bigl(\L[u]\otimes\Z[v]\bigr/(v^2=uv=0)\bigr)^{\otimes m}=R^*(\varDelta^0)^{\otimes m}.
\end{equation}
For $m=1$, the map $s_1\colon E^{0,*}=\k[v]\to E^{-1,*}$ given by
\[
  s_1(a_0+a_1v+\ldots+a_jv^j)=(a_2v+a_3v^2+\ldots+a_jv^{j-1})u
\]
is a cochain homotopy. Indeed, we can write an element of $E$ as
either $x$ or $xu$ with $x=a_0+a_1v+\ldots+a_jv^j\in E^{0,2j}$. In
the former case, $ds_1x=x-a_0-a_1v=x-\iota\varrho x$ and
$s_1dx=0$. In the latter case, $xu\in E^{-1,2j}$, then
$ds_1(xu)=0$ and $s_1d(xu)=xu-a_0u=xu-\iota\varrho(xu)$. In both
cases~\eqref{chaineq1} holds. Now we may assume by induction that
for $m=k-1$ there is a cochain homotopy operator $s_{k-1}\colon
E_{k-1}\to E_{k-1}$. Since $E_k=E_{k-1}\otimes E_1$,
$\varrho_k=\varrho_{k-1}\otimes\varrho_1$ and
$\iota_k=\iota_{k-1}\otimes\iota_1$, a direct check shows that the
map
\[
  s_k=s_{k-1}\otimes\id+\iota_{k-1}\varrho_{k-1}\otimes s_1
\]
is a cochain homotopy between $\operatorname{id}$ and
$\iota_k\varrho_k$, which finishes the proof for
$K=\varDelta^{m-1}$.

In the case of arbitrary $K$ the algebras
$\L[u_1,\ldots,u_m]\otimes\Z[K]$ are $R^*(K)$ are obtained
from~\eqref{emalg} and~\eqref{rmalg} respectively by factoring out
the Stanley--Reisner ideal $\mathcal I_K$. This factorisation does
not affect the properties of the constructed map~$s$, which
finishes the proof.
\end{proof}

Now comparing the additive structure of $R^*(K)$ with that of the
cellular cochains $C^*(K)$, we see that the two coincide:
\begin{lemm}\label{cellcom}
The map
\begin{align*}
  g\colon R^*(K) &\to C^*(\zk),\\
  u_\omega v_\s &\mapsto\mathcal T(\s,\omega)^*
\end{align*}
is an isomorphism of bigraded differential modules. In particular,
we have an additive isomorphism
\[
  H[R^*(K)]\cong H^*(\zk).
\]
\end{lemm}

Having identified the algebra $R^*$ with the cellular cochains of
$\zk$, we can also interpret the cohomology isomorphism from
Lemma~\ref{iscoh} topologically. To do this we shall identify the
Koszul algebra $\Lambda[u_1,\ldots,u_m]\otimes\Z[K]$ with the
cellular cochains of a space homotopy equivalent to $\zk$.

Let $S^\infty$ be an infinite-dimensional sphere obtained as a
direct limit (union) of standardly embedded odd-dimensional
spheres. The space $S^\infty$ is contractible and has a cell
decomposition with one cell in every dimension. The boundary of an
even-dimensional cell is the closure of the appropriate
odd-dimensional cell, while the boundary of an odd cell is zero.
The 2-skeleton of this cell decomposition is a 2-disc decomposed
as shown on Figure~\ref{cellfig}, while the 1-skeleton is the
circle $S^1\subset S^\infty$. The cellular cochain complex of
$S^\infty$ can be identified with the algebra
\[
  \L[u]\otimes\Z[v],\quad\deg u=1,\deg v=2,\quad du=v,dv=0.
\]
From the obvious functorial properties of Construction~\ref{nsc}
we obtain a deformation retraction
\[
  \zk=(D^2,S^1)^K\hookrightarrow(S^\infty,S^1)^K\longrightarrow(D^2,S^1)^K
\]
onto a cellular subcomplex.

The cellular cochains of the $K$-power $(S^\infty,S^1)^K$ can be
identified with the Koszul algebra
$\L[u_1,\ldots,u_m]\otimes\Z[K]$. Since
$\zk\subset(S^\infty,S^1)^K$ is a deformation retract, the
cellular cochain map
\[
  \L[u_1,\ldots,u_m]\otimes\Z[K]=
  C^*\bigl((S^\infty,S^1)^K\bigr)\to C^*(\zk)=R^*(K),
\]
induces an isomorphism in cohomology. In fact, the algebraic
homotopy map $s$ constructed in the proof of Lemma~\ref{iscoh} is
the map induced on the cochains by the topological homotopy.

\subsection{Cellular cochain algebras}
Here we introduce a multiplication for cellular cochains of $\zk$
and establish a ring isomorphism in Lemma~\ref{cellcom}. This task
runs into a complication because cellular cochains in general do
not carry a functorial associative multiplication; the classical
definition of the cohomology multiplication involves a diagonal
map, which is not cellular. However, in our case there is a way to
construct a canonical cellular approximation of the diagonal map
$\D\colon\zk\to\zk\times\zk$ in such a way that the resulting
multiplication in cellular cochains coincides with that
in~$R^*(K)$.

The standard definition of the multiplication in cohomology of a
cell complex $X$ via cellular cochains is as follows. Consider a
composite map of cellular cochain complexes:
\begin{equation}\label{cemul}
\begin{CD}
  C^*(X)\otimes C^*(X) @>\times>> C^*(X\times X)
  @>\widetilde{\D}^*>> C^*(X).
\end{CD}
\end{equation}
Here the map $\times$ assigns to a cellular cochain $c_1\otimes
c_2\in C^{q_1}(X)\otimes C^{q_2}(X)$ the cochain $c_1\times c_2\in
C^{q_1+q_2}(X\times X)$ whose value on a cell $e_1\times e_2\in
X\times X$ is $(-1)^{q_1q_2}c_1(e_1)c_2(e_2)$. The map
$\widetilde{\D}^*$ is induced by a cellular approximation
$\widetilde{\D}$ of the diagonal map $\D\colon X\to X\times X$. In
cohomology, the map~\eqref{cemul} induces a multiplication
$H^*(X)\otimes H^*(X)\to H^*(X)$ which does not depend on a choice
of cellular approximation and is functorial. However, the
map~\eqref{cemul} is not itself functorial because of the
arbitrariness in the choice of a cellular approximation.

In the special case $X=\zk$ we may apply the following
construction. Consider a map $\widetilde{\D}\colon D^2\to
D^2\times D^2$, defined in polar coordinates $z=\rho
e^{i\varphi}\in D^2$, $0\le\rho\le1$, $0\le\varphi<2\pi$ as
follows:
\[
  \rho e^{i\varphi}\mapsto\left\{
  \begin{array}{ll}
    (1+\rho(e^{2i\varphi}-1),1)&\text{ for }0\le\varphi\le\pi,\\
    (1,1+\rho(e^{2i\varphi}-1))&\text{ for }\pi\le\varphi<2\pi.
  \end{array}
  \right.
\]
This is a cellular map taking $\partial D^2$ to $\partial
D^2\times\partial D^2$ and homotopic to the diagonal $\D\colon
D^2\to D^2\times D^2$ in the class of such maps. Taking an
$m$-fold product, we obtain a cellular approximation
\[
  \widetilde{\D}\colon(D^2)^m\to(D^2)^m\times(D^2)^m
\]
which restricts to a cellular approximation for the diagonal
map of $\zk$ for arbitrary $K$, as described by the following
commutative diagram:
\[
\begin{CD}
  \zk @>>> (D^2)^m\\
  @V\widetilde{\D}VV @VV\widetilde{\D}V\\
  \zk\times\zk @>>> (D^2)^m\times(D^2)^m
\end{CD}.
\]

Note that this diagonal approximation is functorial with respect
to those maps $\mathcal Z_{K_1}\to\mathcal Z_{K_2}$ of
moment-angle complexes that are induced by simplicial maps $K_1\to
K_2$ (see Lemma~\ref{mamap}).

\begin{lemm}\label{cellappr}
The cellular cochain algebra $C^*(\zk)$ defined by the diagonal
approximation $\widetilde\D\colon\zk\to\zk\times\zk$
and~\eqref{cemul} is isomorphic to $R^*(K)$. Therefore, we get an
isomorphism of the cohomology algebras:
\[
  H[R^*(K)]\cong H^*(\zk;\Z).
\]
\end{lemm}
\begin{proof}
We first consider the case $K=\D^0$, that is, $\zk=D^2$. The
cellular cochain complex of $D^2$ is additively generated by the
cochains $1\in C^0(D^2)$, $T^*\in C^1(D^2)$ and $D^*\in C^2(D^2)$
dual to the corresponding cells, see Figure~\ref{cellfig}. The
multiplication defined in $C^*(D^2)$ by~\eqref{cemul} is trivial,
so we get a multiplicative isomorphism
\[
  R^*(\D^0)=\L[u]\otimes\Z[v]/(v^2=uv=0)\to C^*(D^2).
\]
Now, for $K=\D^{m-1}$ we obtain a multiplicative isomorphism
\[
  f\colon R^*(\D^{m-1})=
  \L[u_1,\ldots,u_m]\otimes\Z[v_1,\ldots,v_m]/(v_i^2=u_iv_i=0)
  \to C^*((D^2)^m)
\]
by taking the tensor product. Since $\zk\subseteq(D^2)^m$ is a
cell subcomplex for arbitrary $K$ we obtain a multiplicative map
$q\colon C^*((D^2)^m)\to C^*(\zk)$. Now consider the commutative
diagram
\[
\begin{CD}
  R^*(\D^{m-1}) @>f>> C^*((D^2)^m)\\
  @V{p}VV       @VVq V\\
  R^*(K)        @>g>> C^*(\zk).
\end{CD}
\]
Here the maps $p$, $f$ and $q$ are multiplicative, while $g$ is an
additive isomorphism by Lemma~\ref{cellcom}. Take $\alpha,\beta\in
R^*(K)$. Since $p$ is onto, we have $\alpha=p(\alpha')$ and
$\beta=p(\beta')$. Then
\[
  g(\alpha\beta)=gp(\alpha'\beta')=qf(\alpha'\beta')=
  gp(\alpha')gp(\beta')=g(\alpha)g(\beta),
\]
and $g$ is also a multiplicative isomorphism, which finishes the
proof.
\end{proof}

Combining the results of Lemmas \ref{koscom}, \ref{tamap},
\ref{iscoh} and \ref{cellappr}, we come to the main result of this
section.

\begin{theo}\label{zkcoh}
There is an isomorphism, functorial in $K$, of bigraded algebras
\[
  H^{*,*}(\zk;\Z)\cong\Tor_{\Z[v_1,\ldots,v_m]}\bigl(\Z[K],\Z\bigr)\cong
  H\bigl[\L[u_1,\ldots,u_m]\otimes\Z[K],d\bigr],
\]
where the bigrading and the differential in the last algebra are
defined by~\eqref{diff}.
\end{theo}

As an illustration, we give two examples of particular cohomology
calculations, which have a transparent geometrical interpretation.
More examples of calculations may be found in~\cite{bu-pa02}.

\begin{exam}
1. Let $K=\partial\varDelta^{m-1}$. Then
\[
  \Z[K]=\Z[v_1,\ldots,v_m]/(v_1\cdots v_m).
\]
The fundamental class of $\zk\cong S^{2m-1}$ is represented by the
bideg $(-1,2m)$ cocycle $u_1v_2v_3\cdots
v_m\in\Lambda[u_1,\ldots,u_m]\otimes\Z[K]$.

2. Let $K=\{v_1,\ldots,v_m\}$ ($m$ points). Then $\zk$ is homotopy
equivalent to the complement in $\C^m$ to the set of all
codimension-two coordinate planes, see Example~\ref{csaex}. Then
\[
  \Z[K]=\Z[v_1,\ldots,v_m]/(v_iv_j, \ i\ne j).
\]

The subspace of cocycles in $R^*(K)$ is generated by
\[
  v_{i_1}u_{i_2}u_{i_3}\cdots u_{i_k},\quad k\ge2\text{ and }
  i_p\ne i_q\text{ for }p\ne q,
\]
and has dimension $m\binom{m-1}{k-1}$. The subspace of
coboundaries is generated by the elements of the form
\[
  d(u_{i_1}\cdots u_{i_k})
\]
and is $\binom mk$-dimensional. Therefore
\[
\begin{array}{l}
  \dim H^{0}(\zk)=1,\\[2mm]
  \dim H^{1}(\zk)=H^{2}(U(K))=0,\\[2mm]
  \dim H^{k+1}(\zk)=
  m\binom{m-1}{k-1}-\binom mk=(k-1)\binom mk,
  \quad 2\le k\le m,
\end{array}
\]
and
multiplication in
the 
cohomology of $\zk$ is trivial. Note that in general
multiplication in
the 
cohomology of $\zk$ is far from being trivial; for example if $K$ is a
sphere triangulation then $\zk$ is a manifold by Lemma~\ref{zkman}.
\end{exam}

The above cohomology calculation suggests that the complement of
the subspace arrangement from the previous example is homotopy
equivalent to a wedge of spheres. This is indeed the case, as the
following theorem shows.

\begin{theo}[Grbi\'c--Theriault \cite{gr-th04}]
The complement of the set of all codimension-two coordinate
subspaces in $\C^m$ has the homotopy type of the wedge of spheres
\[
  \bigvee_{k=2}^{m} (k-1)\binom{m}{k}S^{k+1}.
\]
\end{theo}
The proof is based on an analysis of the homotopy fibre of the
inclusion $\djs(K)\hookrightarrow BT^m$, which is homotopy
equivalent to $\zk$ (or $U(K)$) by Corollary~\ref{hofib}. We shall
return to coordinate subspace arrangements once again in the next
section.

\section{Applications to combinatorial commutative algebra}

\subsection{
A 
multiplicative version of
Hochster's 
theorem}\label{mhoch}
As a first application we give a proof of a generalisation of
Hochster's theorem (Theorem~\ref{hoch}) obtained by Baskakov
in~\cite{bask02}.

The bigraded structure in
the 
cellular cochains of $\zk$ can be
further refined as
\[
  C^*(\zk)=\bigoplus_{\omega\subseteq[m]}C^{*,\,2\omega}(\zk)
\]
where $C^{*,\,2\omega}(\zk)$ is the subcomplex generated by the
cochains $\mathcal T(\s,\omega\setminus\s)^*$ with
$\s\subseteq\omega$ and $\s\in K$. Thus, $C^*(\zk)$ now becomes a
$\Z\oplus\Z^m$-graded module, and the bigraded cohomology groups
decompose accordingly as
\begin{equation}\label{bomeg}
  H^{-i,\,2j}(\zk)=
  \bigoplus_{\omega\subseteq[m]\colon|\omega|=j}H^{-i,\,2\omega}(\zk)
\end{equation}
where $H^{-i,\,2\omega}(\zk):=H^{-i}[C^{*,\,2\omega}(\zk)]$.

Given two simplicial complexes $K_1$ and $K_2$ with vertex sets
$V_1$ and $V_2$ respectively, their \emph{join} is the following
complex on $V_1\sqcup V_2$:
\[
  K_1*K_2:=\{\s\subseteq V_1\sqcup V_2\colon
  \s=\sigma_1\cup\sigma_2,\:\s_1\in K_1,\:\s_2\in K_2\}.
\]
Now we introduce a multiplication in the sum
\[
  \mathop{\bigoplus_{p\ge-1,}}\limits_{\omega\subseteq[m]}
  \widetilde H^p(K_\omega)
\]
where $K_\omega$ is the full subcomplex and $\widetilde
H^{-1}(\varnothing)=\Z$, as follows. Take two elements
$\alpha\in\widetilde{H}^p(K_{\omega_1})$ and
$\beta\in\widetilde{H}^q(K_{\omega_2})$. Assume that
$\omega_1\cap\omega_2=\varnothing$. Then we have an inclusion of
subcomplexes
\[
  i\colon K_{\omega_1\cup\omega_2}=K_{\omega_1}\sqcup
  K_{\omega_2}\hookrightarrow K_{\omega_1}*K_{\omega_2}
\]
and an isomorphism of reduced simplicial cochains
\[
  f\colon\widetilde{C}^p(K_{\omega_1})\otimes\widetilde{C}^q(K_{\omega_2})
  \stackrel\cong\longrightarrow
  \widetilde{C}^{p+q+1}(K_{\omega_1}*K_{\omega_2}).
\]
Now set
\[
  \alpha\cdot\beta:=\left\{ \begin{array}{ll}
    0, & \omega_1\cap\omega_2\ne\varnothing,\\
    i^*f(a\otimes b)\in\widetilde
    H^{p+q+1}(K_{\omega_1\sqcup\omega_2}), &
    \omega_1\cap\omega_2=\varnothing.
  \end{array}\right.
\]

\begin{theo}[{Baskakov~\cite[Th.~1]{bask02}}]\label{mulho}
There are isomorphisms
\[
  \widetilde H^p(K_\omega)\stackrel\cong\longrightarrow
  H^{p+1-|\omega|,2\omega}(\zk)
\]
which are functorial with respect to simplicial maps and induce a
ring isomorphism
\[
  \gamma\colon
  \mathop{\bigoplus_{p\ge-1,}}\limits_{\omega\subseteq[m]}
  \widetilde H^p(K_\omega)\stackrel\cong\longrightarrow
  H^*(\zk).
\]
\end{theo}
\begin{proof}
Define a map of cochain complexes
\[
  \widetilde C^*(K_\omega)\to C^{*+1-|\omega|,2\omega}(\zk),
  \quad \sigma^*\mapsto\mathcal T(\sigma,\omega\setminus\s)^*.
\]
It is a functorial isomorphism by observation, whence the
isomorphism of the cohomology groups follows.

The statement about the ring isomorphism follows from the
isomorphism $H^*(\zk)\cong H[R^*(K)]$ established in
Lemma~\ref{cellcom} and analysing the ring structure in $R^*(K)$.
\end{proof}

\begin{coro}\label{zkfsc}
There is an isomorphism
\[
  H^{-i,2j}(\zk)\cong\bigoplus_{\omega\subseteq[m]\colon|\omega|=j}
  \widetilde H^{j-i-1}(K_\omega).
\]
\end{coro}
As a further corollary we obtain Hochster's theorem
(Theorem~\ref{hoch}):
\[
  \Tor^{-i,*}_{\Z[v_1,\ldots,v_m]}(\Z[K],\Z)\cong
  \bigoplus_{\omega\subseteq[m]}
  \widetilde H^{|\omega|-i-1}(K_\omega).
\]

\subsection{Alexander duality and coordinate subspace arrangements revisited}
The multiplicative version of Hochster's can also be applied to
cohomology calculations of subspace arrangement complements.

A coordinate subspace can be defined either by setting some
coordinates to zero as in~\eqref{cossp}, or as the linear span of
a subset of the standard basis in $\C^m$. This gives an
alternative way to parametrise coordinate subspace arrangements by
simplicial complexes. Namely, we can write
\[
  \{ L_\omega\colon\omega\notin K \}=
  \{\mathop{\mathrm{span}}\langle e_{i_1},\dots,e_{i_k}\rangle\colon
  \{i_1,\dots,i_k\}\in \widehat{K} \}
\]
where $\widehat{K}$ is the simplicial complex given by
\[
  \widehat{K}:=\{\omega\subseteq[m]\colon[m]\setminus\omega\notin K\}.
\]
It is called the \emph{dual complex} of $K$. The cohomology of
full subcomplexes in $K$ is related to the homology of links in
$\widehat{K}$ by means of the following combinatorial version of
the Alexander duality theorem.

\begin{theo}[Alexander duality]\label{aldua}
Let $K\ne\D^{m-1}$ be a simplicial complex on the set~$[m]$ and
$\s\notin K$, that is, $\widehat\s=[m]\setminus\s\in\widehat K$.
Then there are isomorphisms
\[
  \widetilde{H}_j(K_{\s})\cong
  \widetilde{H}^{|\s|-3-j}(\link_{\widehat{K}}\widehat\s).
\]
In particular, for $\sigma=[m]$ we get
\[
  \widetilde{H}_j(K)\cong
  \widetilde{H}^{m-3-j}(\widehat{K}),\quad-1\le j\le m-2.
\]
\end{theo}
A proof can be found in~\cite[\S2.2]{bu-pa04}. Using the duality
between the full subcomplexes of $K$ and links of $\widehat K$ we
can reformulate the cohomology calculation of~$U(K)$ as follows.

\begin{prop}\label{gmcsa}
We have
\[
  \widetilde{H}_i(U(K))\cong\bigoplus_{\sigma\in\widehat{K}}
  \widetilde{H}^{2m-2|\sigma|-i-2}(\link_{\widehat{K}}\sigma).
\] 
\end{prop}
\begin{proof}
From Proposition~\ref{ukzk} and Corollary~\ref{zkfsc} we obtain
\[
  H_p(U(K))=
  \bigoplus_{\tau\subseteq[m]}\widetilde H_{p-|\tau|-1}(K_\tau).
\]
Nonempty simplices $\tau\in K$ do not contribute to the above sum,
since the corresponding subcomplexes $K_\tau$ are contractible.
Since $\widetilde{H}^{-1}(\varnothing)=\k$ the empty subset of
$[m]$ only contributes $\k$ to $H^0(U(K))$. Hence we may rewrite
the above formula as
\[
  \widetilde{H}_p(U(K))
  =\bigoplus_{\tau\notin K}\widetilde{H}_{p-|\tau|-1}(K_\tau).
\]
Using Theorem \ref{aldua}, we calculate
\[
  \widetilde{H}_{p-|\tau|-1}(K_\tau)=
  \widetilde H^{|\tau|-3-p+|\tau|+1}(\link_{\widehat
  K}\widehat\tau)=
  \widetilde H^{2m-2|\widehat\tau|-p-2}(\link_{\widehat
  K}\widehat\tau),
\] 
where $\widehat\tau=[m]\setminus\tau$ is a simplex in $\widehat
K$, as required.
\end{proof}

Proposition~\ref{gmcsa} is a particular case of the well-known
\emph{Goresky--Macpherson formula}~\cite[Part~III]{go-ma88}, which
calculates the dimensions of the (co)homology groups of an
arbitrary subspace arrangement in terms of its \emph{intersection
poset} (which coincides with the poset of faces of $\widehat K$ in
the case of coordinate arrangements). We see that the study of
moment-angle complexes not only allows us to retrieve the
multiplicative structure of the cohomology of complex coordinate
subspace arrangement complements, but also connects two seemingly
unrelated results, the Goresky--Macpherson formula from the theory
of arrangements and
Hochester's 
formula from combinatorial commutative algebra.

\subsection{Massey products in
the 
cohomology of $\zk$}
Here we address the question of existence of non-trivial Massey
products in the Koszul complex
\[
  [\Lambda[u_1,\ldots,u_m]\otimes\Z[K],d]
\]
of the face ring. Massey products constitute a series of
higher-order operations (or \emph{brackets}) in the cohomology of
a differential graded algebra, with the second-order operation
coinciding with the cohomology multiplication, while the
higher-order brackets are only defined for certain tuples of
cohomology classes. A geometrical approach to constructing
nontrivial triple Massey products in the Koszul complex of the
face ring has been developed by Baskakov in~\cite{bask03} as an
extension of the cohomology calculation in Theorem~\ref{mulho}. It
is well-known that non-trivial higher Massey products obstruct the
\emph{formality} of a differential graded algebra, which in our
case leads to a family on nonformal moment-angle manifolds $\zk$.

Massey products in the cohomology of the Koszul complex of a local
ring $R$ were studied by Golod~\cite{golo62} in connection with
the calculation of the Poincar\'e series of $\Tor_R(\k,\k)$. The
main result of Golod is a calculation of the Poincar\'e series for
the class of rings with vanishing Massey products in the Koszul
complex (including the cohomology multiplication). Such rings were
called \emph{Golod} in~\cite{gu-le69}, where the reader can find a
detailed exposition of Golod's theorem together with several
further applications. The Golod property of face rings was studied
in~\cite{h-r-w99}, where several combinatorial criteria for
Golodness were given.

The difference between our situation and that of Golod is that we
are mainly interested in the cohomology of the Koszul complex for
the face ring of a sphere triangulation~$K$. The corresponding
face ring $\k[K]$ does not qualify for Golodness, as the
corresponding moment-angle complex $\zk$ is a manifold, and
therefore, the cohomology of the Koszul complex of $\k[K]$ must
possess many non-trivial products. Our approach aims at
identifying a class of simplicial complexes with non-trivial
cohomology product but vanishing higher-order Massey operations in
the cohomology of the Koszul complex.

Let $K_i$ be a triangulation of a sphere $S^{n_i-1}$ with
$|V_i|=m_i$ vertices, $i=1,2,3$. Set $m:=m_1+m_2+m_3$,
$n:=n_1+n_2+n_3$, and
\[
\begin{array}{ll}
  K:=K_1*K_2*K_3,\quad
  \zk=\mathcal Z_{K_1}\times\mathcal Z_{K_2}\times\mathcal
  Z_{K_3}.
\end{array}
\]
Note that $K$ is a triangulation of $S^{n-1}$ and $\zk$ is an
$(m+n)$-manifold.

Given $\s\in K$, the \emph{stellar subdivision} of $K$ at $\s$ is
obtained by replacing the star of $\s$ by the cone over its
boundary:
\[
  \zeta_\s(K)=(K\setminus\star_K\s)\cup
  (\cone\partial\star_K\s).
\]

Now choose maximal simplices $\s_1\in K_1$, $\s_2',\s_2''\in K_2$
such that $\s_2'\cap\s_2''=\varnothing$, and $\s_3\in K_3$. Set
\[
  \widetilde K:=\zeta_{\s_1\cup\s_2'}(\zeta_{\s_2''\cup\s_3}(K)).
\]
Then $\widetilde K$ is a triangulation of $S^{n-1}$ with $m+2$
vertices. Take generators
\[
  \beta_i\in\widetilde H^{n_i-1}(\widetilde K_{V_i})
  \cong\widetilde H^{n_i-1}(S^{n_i-1}),\quad i=1,2,3,
\]
where $\widetilde K_{V_i}$ is the restriction of $\widetilde K$ to
the vertex set of $K_i$, and set
\[
  \alpha_i:=\gamma(\beta_i)\in H^{n_i-m_i,2m_i}
  (\mathcal Z_{\widetilde K})\subseteq H^{m_i+n_i}
  (\mathcal Z_{\widetilde K}),
\] 
where $\gamma$ is the isomorphism from Theorem~\ref{mulho}. Then
\[
  \beta_1\beta_2\in\widetilde H^{n_1+n_2-1}
  (\widetilde K_{V_1\sqcup V_2})\cong\widetilde
  H^{n_1+n_2-1}(S^{n_1+n_2-1}\setminus\text{pt})=0,
\]
and therefore, $\alpha_1\alpha_2=\gamma(\beta_1\beta_2)=0$, and
similarly $\alpha_2\alpha_3=0$. In these circumstances the triple
Massey product $\langle \alpha_1,\alpha_2,\alpha_3\rangle\subset
H^{m+n-1}(\mathcal Z_{\widetilde K})$ is defined. Recall that
$\langle\alpha_1,\alpha_2,\alpha_3\rangle$ is the set of
cohomology classes represented by the cocycles $(-1)^{\deg
a_1+1}a_1f+ea_3$ where $a_i$ is a cocycle representing $\alpha_i$,
$i=1,2,3$, while $e$ and $f$ are cochains satisfying $de=a_1a_2$,
$df=a_2a_3$. A Massey product is called \emph{trivial} if it
contains zero.

\begin{theo}\label{ntmas}
The triple Massey product
\[
  \langle \alpha_1,\alpha_2,\alpha_3\rangle\subset
  H^{m+n-1}(\mathcal Z_{\widetilde K})
\]
in the cohomology of $(m+n+2)$-manifold $\mathcal Z_{\widetilde
K}$ is non-trivial.
\end{theo}
\begin{proof}
Consider the subcomplex of $\widetilde K$ consisting of those two
new vertices added to $K$ in the process of stellar subdivision.
By Lemma~\ref{mamap}, the inclusion of this subcomplex induces an
embedding of a 3-dimensional sphere $S^3\subset\mathcal
Z_{\widetilde K}$. Since the two new vertices are not joined by an
edge in $\mathcal Z_{\widetilde K}$, the embedded 3-sphere defines
a non-trivial class in $H^3(\mathcal Z_{\widetilde K})$. By
construction the dual cohomology class is contained in the Massey
product $\langle\alpha_1,\alpha_2,\alpha_3\rangle$. On the other
hand, this Massey product is defined up to elements from the
subspace
\[
  \alpha_1\cdot H^{m_2+m_3+n_2+n_3-1}(\mathcal Z_{\widetilde K})+
  \alpha_3\cdot H^{m_1+m_2+n_1+n_2-1}(\mathcal Z_{\widetilde K}).
\]
The multigraded components of the group
$H^{m_2+m_3+n_2+n_3-1}(\mathcal Z_{\widetilde K})$ 
different from that determined by the full subcomplex $\widetilde
K_{V_2\sqcup V_3}$ do not affect the nontriviality of the Massey
product, while the multigraded component corresponding to
$\widetilde K_{V_2\sqcup V_3}$ is zero since this subcomplex is
contractible. The group $H^{m_1+m_2+n_1+n_2-1}(\mathcal
Z_{\widetilde K})$ is treated similarly. It follows that the
Massey product contains a unique nonzero element in its
multigraded component and so is nontrivial.
\end{proof}

As is well known, the nontriviality of Massey products obstructs
formality of manifolds, see e.g.~\cite{ba-ta00}.

\begin{coro}
For every sphere triangulation $\widetilde K$ obtained from
another triangulation by applying two stellar subdivisions as
described above, the 2-connected moment-angle manifold $\mathcal
Z_{\widetilde K}$ is nonformal.
\end{coro}

In the proof of Theorem~\ref{ntmas} the nontriviality of the
Massey product is established geometrically. A parallel argument
may be carried out algebraically in terms of the algebra $R^*(K)$,
as illustrated in the following example.

\begin{exam}
Consider the simple polytope $P^3$ shown on
Figure~\ref{masseyfig}.
\begin{figure}[h]
\begin{center}
\begin{picture}(75,80)
  \put(10,15){\line(2,-1){10}}
  \put(10,15){\line(0,1){40}}
  \put(20,50){\line(-2,1){10}}
  \put(20,50){\line(0,-1){40}}
  \put(20,50){\line(1,0){20}}
  \put(40,50){\line(1,-1){5}}
  \put(45,45){\line(0,-1){35}}
  \put(45,10){\line(-1,0){25}}
  \put(45,10){\line(1,1){20}}
  \put(65,30){\line(0,1){35}}
  \put(65,65){\line(-1,-1){20}}
  \put(40,50){\line(1,1){20}}
  \put(60,70){\line(1,-1){5}}
  \put(60,70){\line(-1,0){35}}
  \put(25,70){\line(-1,-1){15}}
  \multiput(25,30)(5.1,0){8}{\line(1,0){3}}
  \multiput(25,30)(0,5.1){8}{\line(0,1){3}}
  \multiput(25,30)(-5.2,-5.2){3}{\line(-1,-1){3.6}}
  \put(14,30){$w_1$}
  \put(49,55){$w_2$}
  \put(30,58){$v_1$}
  \put(35,5){\vector(0,1){4}}
  \put(34,3){$v_2$}
  \put(54,35){$v_3$}
  \put(13,66){\vector(1,-1){4}}
  \put(10,67){$v_4$}
  \put(30,32){$v_5$}
  \put(50,75){\vector(0,-1){4}}
  \put(48,76){$v_6$}
\end{picture}
\caption{ } \label{masseyfig}
\end{center}
\end{figure}
This polytope is obtained by cutting two non-adjacent edges off a
cube and has 8 facets. The dual triangulation $K_P$ is obtained
from an octahedron by applying stellar subdivisions at two
non-adjacent edges. The face ring is
\[
  \Z[K_P]=\Z[v_1,\ldots,v_6,w_1,w_2]/\mathcal I_{K_P},
\]
where $v_i$, $i=1,\ldots,6$, are the generators coming from the
facets of the cube and $w_1,w_2$ are the generators corresponding
to the two new facets, see Figure~\ref{masseyfig}, and
\[
  \mathcal I_P=(v_1v_2,v_3v_4,v_5v_6,w_1w_2,v_1v_3,v_4v_5,
  w_1v_3,w_1v_6,w_2v_2,w_2v_4).
\]
The corresponding algebra $R^*(K_P)$ has additional generators
$u_1,\ldots,u_6,t_1,t_2$ of total degree 1 satisfying $du_i=v_i$
and $dt_i=w_i$. Consider the cocycles
\[
  a=v_1u_2,\quad b=v_3u_4,\quad c=v_5u_6
\]
and the corresponding cohomology classes $\alpha,\beta,\gamma\in
H^{-1,4}[R^*(K)]$. The equations
\[
  ab=de,\quad bc=df
\]
have a solution $e=0$, $f=v_5u_3u_4u_6$, so the triple Massey
product $\langle\alpha,\beta,\gamma\rangle\in H^{-4,12}[R^*(K)]$
is defined. This Massey product is nontrivial by
Theorem~\ref{ntmas}. The cocycle
\[
  af+ec=v_1v_5u_2u_3u_4u_6
\]
represents a nontrivial cohomology class
$[v_1v_5u_2u_3u_4u_6]\in\langle\alpha,\beta,\gamma\rangle$ and so
the algebra $R^*(K_P)$ and the manifold $\mathcal Z_{K_P}$ are not
formal.
\end{exam}

In view of Theorem~\ref{ntmas},
the question arises 
of describing the class of simplicial complexes $K$ for which the
algebra $R^*(K)$ (equivalently, the Koszul algebra
$[\Lambda[u_1,\ldots,u_m]\otimes\Z[K],d]$ or the space $\zk$) is
formal (in particular, does not contain nontrivial Massey
products). For example, a direct calculation shows that this is
the case if $K$ is the boundary of a polygon.

\subsection{Toral rank conjecture}
Here we relate our cohomological calculations with moment-angle
complexes to an interesting conjecture in the theory of
transformation groups. This `toral rank conjecture' has strong
links with rational homotopy theory, as described
in~\cite{al-pu93}. Therefore this last subsection, although not
containing new results, aims at encouraging rational homotopy
theorists to turn their attention to combinatorial commutative
algebra of simplicial complexes.

A torus action on a space $X$ is called \emph{almost free} if all
isotropy subgroups are finite. The \emph{toral rank} of $X$,
denoted $\trk(X)$, is the largest $k$ for which there exists an
almost free $T^k$-action on $X$.

The \emph{toral rank conjecture} of Halperin~\cite{halp85}
suggests that
\[
  \dim H^*(X;\Q)\ge2^{\trk(X)}
\]
for any finite dimensional space $X$. Equality is achieved, for
example, if $X=T^k$.

Moment angle complexes provide a wide class of almost free torus
actions:

\begin{theo}[{Davis--Januszkiewicz~\cite[7.1]{da-ja91}}]\label{trank}
Let $K$ be an $(n-1)$-di\-men\-si\-o\-nal simplicial complex with
$m$ vertices. Then $\trk\zk\ge m-n$.
\end{theo}
\begin{proof}
Choose an lsop in $t_1,\ldots,t_n$ in $\Q[K]$ according to
Lemma~\ref{noether} and write
\[
  t_i=\lambda_{i1}v_1+\ldots+\lambda_{im}v_m,\quad i=1,\ldots,n.
\]
Then the matrix $\Lambda=(\lambda_{ij})$ defines a linear map
$\lambda\colon\Q^m\to\Q^n$. Changing $\lambda$ to $k\lambda$ for a
sufficiently large $k$ if necessary, we may assume that $\lambda$
is induced by a map $\Z^m\to\Z^n$, which we also denote
by~$\lambda$. It follows from Lemma~\ref{restr} that for every
simplex $\sigma\in K$ the restriction
$\lambda|_{\Z^\s}\colon\Z^\s\to\Z^n$ of the map $\lambda$ to the
coordinate subspace $\Z^\s\subseteq\Z^m$ is injective.

Denote by $T_\Lambda$ the subgroup in $T^m$ corresponding to the
kernel of the map $\lambda\colon\Z^m\to\Z^n$. Then $T_\Lambda$ is
a product of an $(m-n)$-dimensional torus $N$ and a finite group.
The intersection of the torus $N$ with the coordinate subgroup
$T^{\sigma}\subseteq T^m$ is a finite subgroup. Since the isotropy
subgroups of the $T^m$-action on $\zk$ are of the form $T^\s$ (see
the proof of Theorem~\ref{zkquo}), the torus $N$ acts on $\zk$
almost freely.
\end{proof}

Note that by construction the space $\zk$ is 2-connected.

In view of Theorem~\ref{mulho}, we get the following reformulation
of the toral rank conjecture for $\zk$:
\[
  \dim\bigoplus_{\omega\subseteq[m]}\widetilde
  H^*(K_\omega;\Q)\ge2^{m-n}
\]
for any simplicial complex $K^{n-1}$ on $m$ vertices.

\begin{exam}
Let $K$ the boundary of an $m$-gon. Then the calculation
of~\cite[Exam.~7.22]{bu-pa02} shows that
\[
  \dim H^*(\zk)=(m-4)2^{m-2}+4\ge2^{m-2}.
\]
\end{exam}

\end{document}